\documentclass[12pt,aps,pre,preprint]{revtex4}%

\usepackage{amssymb}
\usepackage{amsmath}
\usepackage[dvips]{graphicx}
\usepackage{colordvi}
\usepackage{times}

\newcommand{\bea}{\begin{eqnarray}}
\newcommand{\eea}{\end{eqnarray}}
\newcommand{\be}{\begin{equation}}
\newcommand{\ee}{\end{equation}}
\newcommand{\un}[1]{\underline{#1}}
\newcommand{\zL}{\Lambda}
\newcommand{\zd}{\delta}

\begin{document}

\title{Equilibrium, fluctuation relations and transport for irreversible deterministic dynamics}

\author{Matteo Colangeli \email{colangeli@calvino.polito.it}}
\affiliation{Dipartimento di Matematica, Politecnico di Torino, Corso Duca degli Abruzzi
24, I-10129 Torino, Italy}

\author{Lamberto Rondoni}
\affiliation{Dipartimento di Matematica, Politecnico di Torino, Corso Duca degli Abruzzi
24, I-10129 Torino, Italy.\\
INFN, Sezione di Torino, Via P. Giura 1, I-10125, Torino, Italy}


\begin{abstract}
In a recent paper [M. Colangeli \textit{et al.}, J.\ Stat.\ Mech.\ P04021, (2011)] it was argued that the
Fluctuation Relation for the phase space contraction rate $\Lambda$ could suitably be extended to
non-reversible dissipative systems. We strengthen here those arguments, providing analytical and numerical
evidence based on the properties of a simple irreversible nonequilibrium baker model. We also consider
the problem of response, showing that the transport coefficients are not affected by the irreversibility
of the microscopic dynamics. In addition, we prove that a form of \textit{detailed balance}, hence of equilibrium,
holds in the space of relevant variables, despite the irreversibility of the phase space dynamics.
This corroborates the idea that the same stochastic description, which arises from a projection onto a
subspace of relevant coordinates, is compatible with quite different underlying deterministic dynamics.
In other words, the details of the microscopic dynamics are largely irrelevant, for what concerns
properties such as those concerning the Fluctuation Relations, the equilibrium behaviour and the
response to perturbations.
\end{abstract}

\maketitle

\section{Introduction}

Physical laws describing the evolution of the microscopic constituents of macroscopic objects
are commonly assumed to be time reversal invariant, and as such they are usually considered in statistical
mechanics. We call T-symmetric an evolution which is time reversal invariant. Differently, the laws of
Thermodynamics, which describe the macroscopic realm, are irreversible.
The theoretical investigation of the so-called Fluctuation Relations (FRs), which began in the early 1990's,
produced one quantitative approach to emergence of irreversible behavior from the statistics of dissipative,
but time reversal invariant, dynamics \cite{Gall,SRE}. The FRs have then been extended and applied to
numerous nonequilibrium phenomena \cite{Bett,RonMejia,Gawedzky,Aurell}. The usual derivation of the FRs
for deterministic dynamics relies on the time reversal invariance of the dynamical equations; likewise, the
derivation of FRs for stochastic dynamics rests on a form of ``reversibility'' which is relevant to the
mesoscopic level of description. Stochastic dynamics, thought to represent a projection of the phase space
dynamics on the space of physical observables (e.g.\ the 1-particle space of a many particle system)
is intrinsically irreversible but, under certain conditions, it is characterized by fluctuations which
allow every path in the state space, which visits certain states in a given chronological order, to be
coupled with with another path, which visits the same states in reverse order. This property can the be
used in the derivation of the FRs for stochastic dynamics.
More precisely, the T-symmetry present in the stochastic approach to the FRs takes the form of a local
detailed balance condition, for nonequilibrium systems driven to a steady state, \cite{CMW,LeboSpohn}.

Recently, the possibility of relaxing the T-symmetry, which is normally assumed to concern the whole phase
space or the whole state space, have been investigated. In Ref.\cite{Gall98} time reversal has been
considered in cases in which the dynamics evolves on a lower dimensional manifold, which is not time
reversal invariant, but is embedded in the phase space which is. In Refs.\cite{MaesIrr,Jona} a similar
investigation is performed for stochastic evolutions, in order to identify the minimal
ingredients required to obtain the FRs or relations of irreversible thermodynamics, such
as the Onsager-Casimir Relations, previously based on time reversal invariance.

In a previous paper, \cite{CKDR}, we too considered the question of how fundamental the T-symmetry
is, in general. To that purpose, we investigated the effect of a ``tunable'' source of irreversibility on
a variation of the deterministic and T-reversible toy model known as dissipative Baker Map \cite{Dorfman}.
In this paper, extending the work of \cite{CKDR}, we test the robustness of the FR and of the linear regime,
with respect to perturbations of the reversibility of the ``microscopic'' dynamics. This uses some of the
insights of the previous work \cite{CKDR}, which showed that smoothness of the invariant measure and the
continuity of the time-reversal symmetry are not necessary in order to derive a FR.
The conjecture which our data are meant to support, from the point of view of a very idealized setting,
is that the time evolution of macroscopic observables of real physical systems is compatible with many
possible underlying dynamics, including irreversible ones. This depends of course on the kinds of
observables at hand.
Hence, we adopted a very simple model in which, however, one may distinguish \textit{relevant} from
\textit{irrelevant} observables, as far as the phase space contraction, denoted by $\zL$, is concerned.
Then, we introduced a source of irreversibility in the dynamics, which affects only the irrelevant variable
but which affects the structure of the steady state probability distributions, including the ``equilibrium''
ones, and we asked what consequences that may have on the validity of the $\zL$-FR, namely the FR for the
quantity $\zL$.\\
In our model, the T-symmetry is relaxed by considering a dynamics which is non-invertible. Mappings which are not homeomorphisms constitute a subclass of non-reversible mappings. We expect our results to hold, more generally, also in some non-reversible but invertible dynamics. To establish a link between deterministic and stochastic dynamical systems, it proves useful to treat systems which are non-reversible in the usual deterministic sense (to be specified below), but which may be endowed with a weaker notion of reversibility (also to be specified below). For example, one may consider a suitable weak notion of reversibility in a projected dynamics (cf. Sec. \ref{sec:sec2}). In particular, our low dimensional toy model may be useful to illustrate how such a weakly reversible projected dynamics is consistent with many, possibly irreversible, deterministic dynamics.\\
Our paper is organized as follows.\\
In Sec.\ref{sec:sec1}, we describe the general setting and focus on the equilibrium version of the model.
The proposed distinction between ``relevant" and ``irrelevant" degrees of freedom, with respect to the
specific observable $\zL$, enable us to verify that the equilibrium version of the
$\zL$-FR is not affected by our source of irreversibility.
In Sec.\ref{sec:sec2}, we derive a detailed balance condition in the projected ``relevant'' space, the one
that would correspond to the $\mu$-space of kinetic theory, in the case of a system of particles. This notion
of detailed balance stems from the more general condition of equilibrium in the full phase space, and is proven
to work even in presence of irreversible dynamics.
In Sec.\ref{sec:sec3}, we consider a nonequilibrium version of our irreversible deterministic map, and study
the validity of the $\zL$-FR, of linear response and of Green-Kubo-like relations for such a map, comparing
ours with other approaches. Conclusions are drawn in Sec. \ref{sec:sec4}.

Our results can be summarized as follows:
\begin{itemize}
\item The condition of \textit{detailed balance}, in the projected relevant space (e.g.\ the
$\mu$-space in kinetic theory), may be derived even from an irreversible equilibrium deterministic
dynamics in phase space. This shows that a stochastic process derived from a suitable projection
onto a proper subspace is compatible with many different underlying deterministic dynamics.
\item Analytical and numerical results prove the validity of the $\Lambda$-FR for a deterministic
dynamical system which is not time reversal invariant and, yet, satisfies a milder, stochastic-like,
notion of reversibility. This form of reversibility merely requires the existence of pairs of
conjugated paths in phase space giving rise to opposite phase space contractions.
\item Linear response theory does not hold in our simple models, which violate the conditions
required by the methods discussed in Refs.\cite{SRE,Gall2} and, in particular, do not enjoy any
property of the kind of Local Thermodynamic Equilibrium. However, the transport coefficients for
the corresponding reversible and irreversible dynamics coincide.
\end{itemize}

\section{Global conservativity vs. local dissipativity}
\label{sec:sec1}
\begin{figure}
   \centering
   \includegraphics[width=11.5cm]{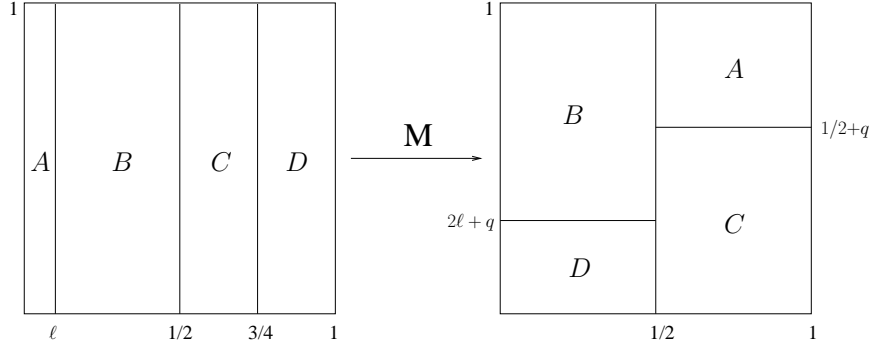}
   \caption{The map $M$ in Eq. (\ref{Map}) for general values of the parameters $\ell$ and $q$.}\label{M4}
\end{figure}

Let us introduce the dynamical system $(\mathcal{U},M,\mu)$, with phase space
$\mathcal{U}:=\mathbb{T}^2:=\mathbb{R}^2/\mathbb{Z}^2$ and mapping $M: \mathcal{U} \to \mathcal{U}$
defined by:
\be
\left(\begin{array}{c}
  x_{n+1} \\
  y_{n+1}
\end{array}\right)
  = M \cdot
\left(\begin{array}{c}
  x_{n} \\
  y_{n}
  \end{array}\right)=\left\{
  \begin{array}{l c}
    \left(\begin{array}{c}
            \dfrac{1}{2\ell}x_{n}+\dfrac{1}{2}\\ \\
            \left(\dfrac{1}{2}-q\right) y_{n} +\left(\dfrac{1}{2}+q\right)
          \end{array}\right) & \quad \text{for $0\leq x_{n} < \ell$}\\
 & \\
        \left(\begin{array}{c}
            \dfrac{x_{n}}{1-2\ell}-\dfrac{\ell}{1-2\ell} \\ \\
            (1-2\ell-q) y_{n} + (2\ell+q)
          \end{array}\right) & \quad \text{for $\ell\leq x_{n} < \frac{1}{2}$} \\
 & \\
        \left(\begin{array}{c}
            2x_{n}-\dfrac{1}{2} \\ \\
            \left(\dfrac{1}{2}+q\right) y_{n}
          \end{array}\right) & \quad \text{for $\frac{1}{2}\leq x_{n} < \frac{3}{4}$} \\
 & \\
        \left(\begin{array}{c}
            2x_{n}-\dfrac{3}{2} \\ \\
            (2\ell+q) y_{n}
          \end{array}\right) & \quad \text{for $\frac{3}{4}\leq x_{n} \leq 1$}
  \end{array}\right. \label{Map}
\ee
for $q \in[0,\frac{1}{2}]$ and $\ell\in(0, \frac{1}{4}]$, cf. Fig. \ref{M4}, and with natural measure $\mu$.

The Jacobian determinant of this map takes the values:
\be
J_{M}(\underline{x})=\left\{\begin{array}{l c}
                     J_A= \dfrac{1}{4\ell}-\dfrac{q}{2\ell} & \quad \text{for $0\leq x < \ell$}\\
 & \\
                     J_B= 1-\dfrac{q}{1-2\ell} & \quad \text{for $\ell\leq x < \dfrac{1}{2}$}\\
 & \\
                     J_C= 1+2q& \quad \text{for $\dfrac{1}{2}\leq x < \dfrac{3}{4}$}\\
 & \\
                     J_D= 4\ell+2q & \quad \text{for $\frac{3}{4}\leq x \leq 1$} \label{J4}
                   \end{array}\right.
\ee
in the four different regions of $\mathcal{U}$.
This model generalizes the one introduced in \cite{CKDR}, as it features the two parameters $\ell$ and $q$,
which can be tuned to produce different forms of ``equilibrium'', i.e. of natural measures $\mu$, which
are called non-dissipative steady states because are characterized by vanishing phase space contraction rates.
In our case, it proves convenient to determine the projection of the invariant probability density on the
$x$-coordinate of this map. This can be accomplished by integrating over the $y$-direction the Perron-Frobenius equation \cite{Lasota} for the measures defined on the square. This, indeed, yields the evolution equation
for the probability measures defined on the $x$ axis which are evolved by the map of the interval $[0,1]$
obtained by projecting $M$ on the $x$ axis.

\begin{figure}
   \centering
   \includegraphics[width=9.5cm]{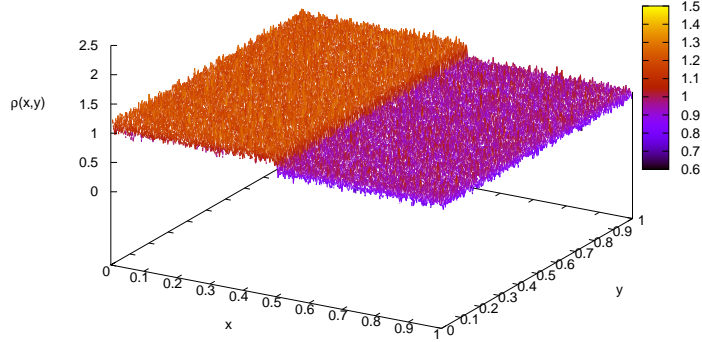}
   \caption{Result of a numerical simulation for the invariant density of the equilibrium reversible map derived from Eq. (\ref{Map}) with $q=0$ and $\ell=0.15$, obtained by evolving a set of $2\cdot10^7$ initial conditions randomly (and uniformly) chosen on the square $[0,1]\times [0,1]$. This shows an invariant density $\rho(x,y)$ which is uniform along the $y$-coordinate and piecewise constant along the $x$-coordinate, attaining the values $\rho_l(x)$ and $\rho_r(x)$ for, respectively $x \in [0,1/2)$ and $x\in[1/2,1]$, given in Eq.(\ref{SRB}).}\label{HistB4rev}
\end{figure}

The calculation of the marginal invariant probability measure can the be performed by
introducing a Markov partition of the unit square, consisting of two regions:
$[0,1/2)$ and $[1/2,1]$, respectively furnished with the invariant densities
$\rho_{l}(x)$ and $\rho_{r}(x)$. Then, the transfer operator $T$ associated
with the projected dynamics, via the projected Perron-Frobenius equation, can be written as:
\be
\left(\begin{array}{c}
        \rho_{l}(x_{n+1}) \\
        \rho_{r}(x_{n+1})
      \end{array}\right)=T\cdot \left(\begin{array}{c}
        \rho_{l}(x_{n}) \\
        \rho_{r}(x_{n})
      \end{array}\right) \label{PF}
\ee
where $T$ is defined by:
\be
T=\left(
    \begin{array}{cc}
      1-2\ell & 1/2 \\
      2\ell & 1/2 \\
    \end{array}
  \right) \label{toptm}
\ee
The matrix $T$ satisfies the Perron-Frobenius Theorem, hence its largest eigenvalue, $\lambda=1$, is
separated from a spectral gap from its other eigenvalue.
Then, the calculation proceeds by evaluating the eigenvectors of the transfer operator corresponding to the dominant eigenvalue. The result of this procedure shows that the invariant probability density of the map (\ref{Map}), projected onto the $x$-axis, depends on the value of $\ell$, but not on $q$, because $q$ only affects
the dynamics along the vertical direction. The corresponding projected density, cf. Fig.\ref{HistB4rev}, is given by the piecewise
constant function:
\be
\rho(x)=\left\{\begin{array}{cc}
                 \rho_l(x)=\dfrac{2}{1+4\ell} & \quad \text{for $0\leq x < \dfrac{1}{2}$} \\ \\
                 \rho_r(x)=\dfrac{8\ell}{1+4\ell} & \quad \text{for $\dfrac{1}{2}\leq x \leq 1$}
               \end{array}\right.
 \; ,\;  \label{SRB}
\ee

The marginal probability density $\rho$ suffices to compute the statistical properties of phase
functions such as the phase space contraction rate $\zL(x,y)=-\log J(x,y)$, because of the special
form of the Jabobian determinants (\ref{J4}), which depend on the $x$-coordinate only. In this case,
one has:
\be
\langle \Lambda \rangle =-\int_{\mathcal{U}} \log J(x,y)\mu(dx\times dy)=
-\int_0^1 \log J(x)\rho(x)dx \nonumber
\ee
\begin{figure}
\begin{center}
  \includegraphics[width=0.80\textwidth]{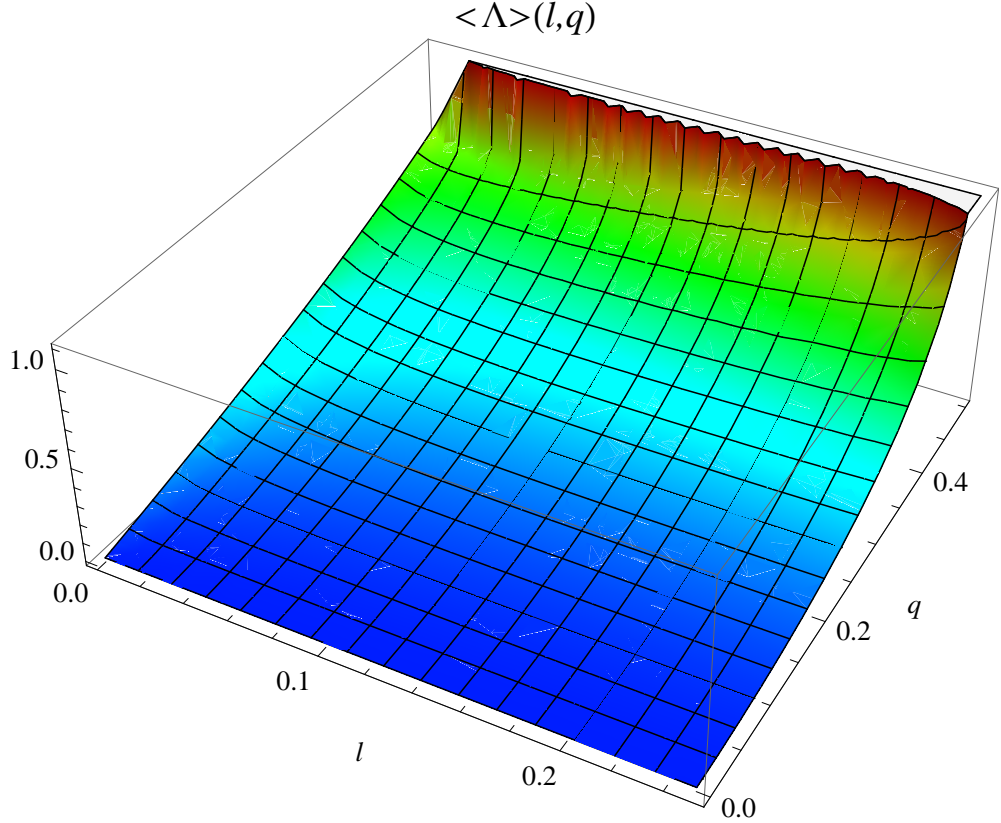}\\
  \caption{The average phase space contraction rate $\langle \Lambda \rangle$, as a function of the parameters
  $\ell$ and $q$.}\label{Lambda}
  \end{center}
\end{figure}
The quantity $\langle \Lambda \rangle$ is represented in Fig.\ref{Lambda} as a function of the parameters
$(\ell,q)$. The figure shows that ``equilibrium'', i.e.\ by definition the condition in which
$\langle \Lambda \rangle$ vanishes, holds only for $q=0$, independently of the value of $\ell$.
The dynamical system of Ref.\cite{CKDR} can be seen as a special case of our map, in which
$q=\frac{1}{2}-2\ell$ and, correspondingly, the only possible equilibrium state of that map is given
by the further choice $\ell=\frac{1}{4}$.

In this Section, we focus on the case $q=0$ and begin by considering the average phase space contraction rate which, in this case,
vanishes $\forall \ell\in(0,\frac{1}{4}]$ and can be written as:
\be
\langle \Lambda \rangle =-\ln J_A \rho_l \ell-\ln J_D \frac{\rho_r}{4}=0 \label{oursigma}
\ee

\begin{figure}
  \centering
   \includegraphics[width=10.5cm]{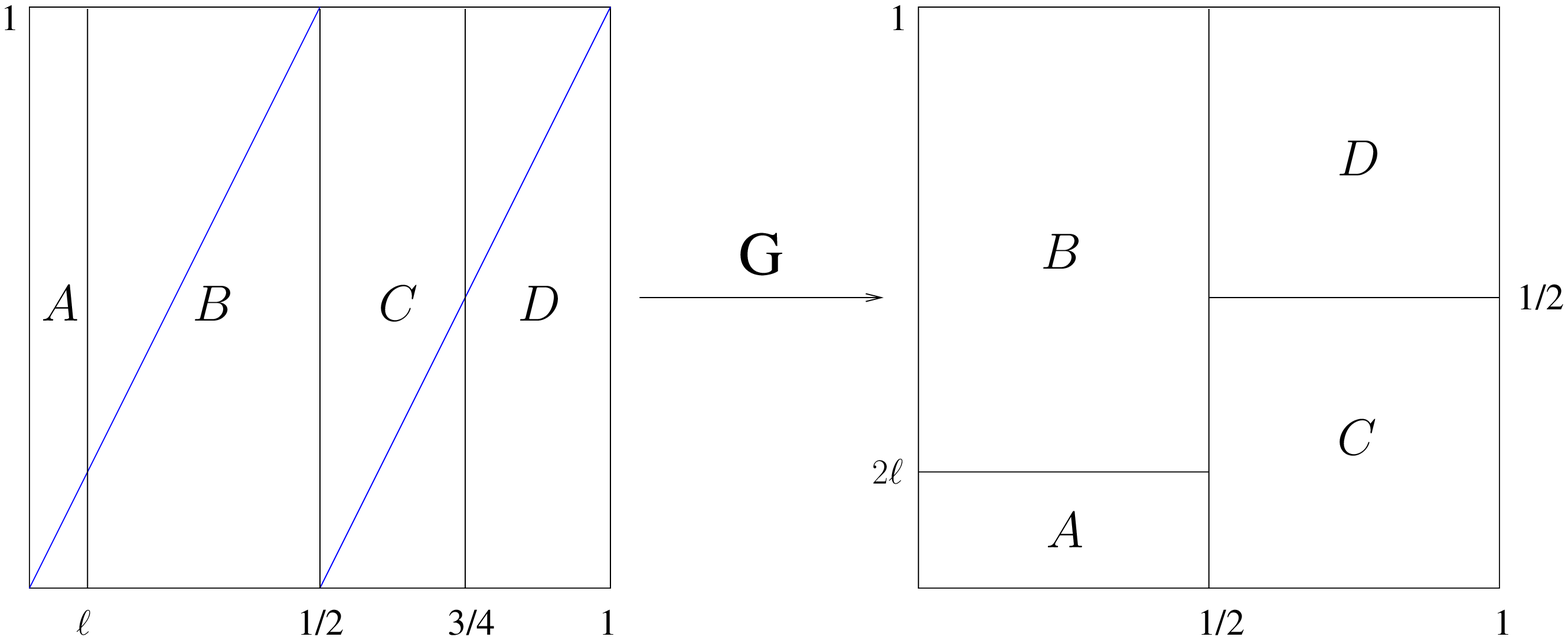}
   \caption{Involution $G$ defined in Eq. (\ref{G}). \textit{Blue lines}: the two diagonals, along which the map $G$ reflects the two halves of the phase space.}\label{involM4}
\end{figure}

The choice $\ell=\frac{1}{4}$, in particular, ensures that the mapping is locally conservative, i.e.\
that $\Lambda(x,y)=0$ uniformly on $\mathcal{U}$, as all Jacobians (\ref{J4}) are unitary.
In particular, $\ell=\frac{1}{4}$ leads to a uniform invariant measure $\mu$, which we call
``microcanonical'' by analogy with the statistical mechanics of an isolated particle system
with given total energy. The projection of $\mu$ is then uniform along the $x$-axis.

For $\ell \neq \frac{1}{4}$ and $q=0$, the invariant density along the $x$-axis remains smooth,
except at one point of discontinuity, $x=1/2$ and, in spite of the fluctuations of the phase
space volumes, the invariant measure is still uniform along the stable manifold,
i.e. the vertical direction, as illustrated by a numerical simulation reported in Fig.\ref{HistB4rev}.
Thus, for $\ell \neq \frac{1}{4}$, we obtain a form of \textit{fluctuating equilibrium}, which we
call ``canonical'' by analogy with the statistical mechanics of a particle system in equilibrium
with a thermostat at a given temperature \cite{JR}.

Let us also observe that our equilibrium dynamics ($q=0$) are time reversal invariant, according
to the standard dynamical systems notion of reversibility \cite{RQ}, because there exists an
involution $G:\mathcal{U}\rightarrow \mathcal{U}$ such that
\be
M G M= G \label{general}
\ee
which attains the form:

\be
\left(\begin{array}{c}
  x_{G} \\
  y_{G}
\end{array}\right)
  = G \cdot
\left(\begin{array}{c}
  x \\
  y
  \end{array}\right)=\left\{
  \begin{array}{l c}
    \left(\begin{array}{c}
            2 x \\ \\
            \dfrac{1}{2} y
          \end{array}\right) & \quad \text{for $0\leq x < \frac{1}{2}$}\\
 & \\
        \left(\begin{array}{c}
            2 x-1 \\ \\
            \frac{1}{2}(y +1)
          \end{array}\right) & \quad \text{for $\frac{1}{2}\leq x \leq 1$}
  \end{array}\right. \label{G}
\ee

\begin{figure}
\begin{center}
  \includegraphics[width=\textwidth]{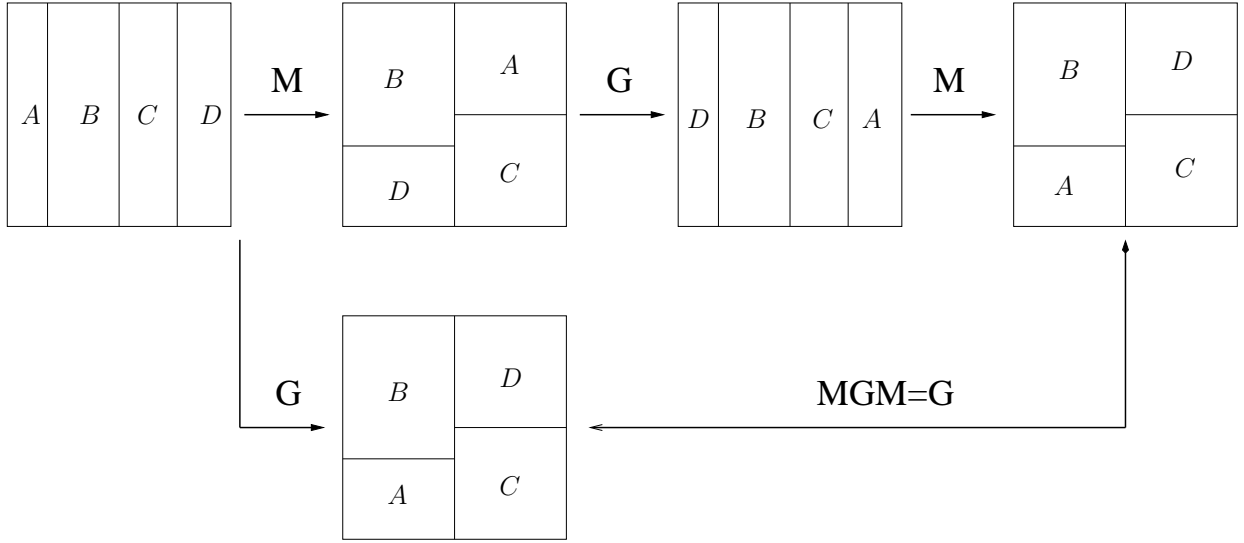}\\
  \caption{Check of reversibility for the map (\ref{Map}).}\label{checkinvolM4}
    \end{center}
\end{figure}

The mapping $G$ in Eq. (\ref{G}) reflects the half squares $[0,1/2)$ and $[1/2,1]$
along the respective diagonals, drawn from their lower left to their upper right corners,
cf. Fig. \ref{involM4}.\\
Hence, according to the definition Eq.(\ref{general}), the baker model in (\ref{Map}) with $q=0$
is T-symmetric for all values of the parameter $\ell$, as shown in Fig. \ref{checkinvolM4}.

\begin{figure}
  \begin{center} \includegraphics[width=12cm,
  width=0.90\textwidth]{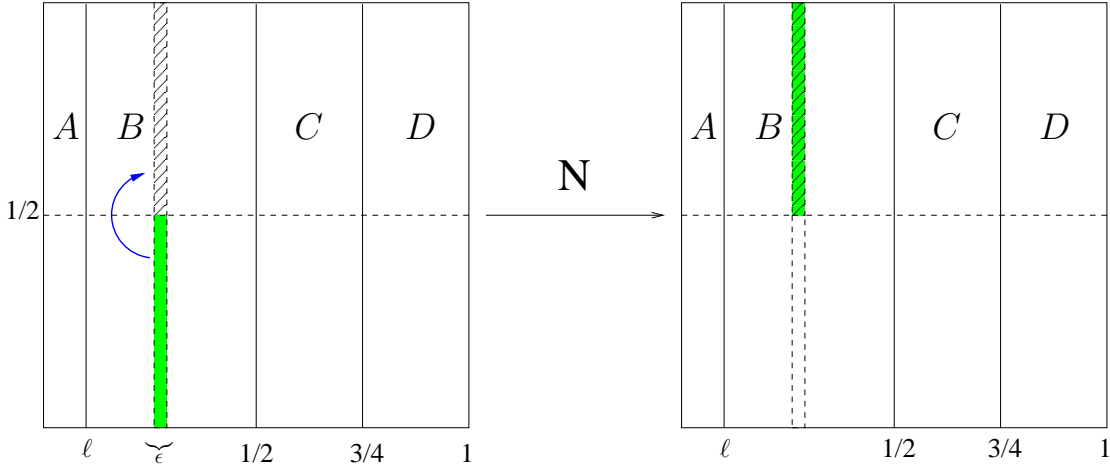}\\
\caption{The map $N$ defined in  Eq.(\ref{N}), which spoils the reversibility of the model}\label{correction} \end{center}
\end{figure}

Consider, now, a trajectory of $n$ time steps,
$\{\un{x}_0,\un{x}_1,...,\un{x}_n\}$, along which the average phase space contraction rate
is given by:
\be
\overline{\zL}_{n}(\un{x}_{0})= - \frac{1}{n}\sum_{k=0}^{n-1}\ln J_M (M^{k}\un{x}_{0}) \label{forward}
\ee
The average phase space contraction rate over the time reversed path is given by:
\bea
\overline{\zL}_{n}(G M^{n}(\un{x}_{0}))&=& - \frac{1}{n}\sum_{k=0}^{n-1}\ln J_M (M^{k}G M^{n}\un{x}_{0})\nonumber\\
&=&-\frac{1}{n}\sum_{k=0}^{n-1}\ln J_M (G M M^{(n-1)-k}\un{x}_{0})\label{backward}
\eea
Then, the relation $J_{M}(\un{x})=J_M^{-1}(G M \un{x})$, yields the known result \cite{CKDR}:
$$
\overline{\zL}_{n}(G M^{n}\un{x}_{0})=-\overline{\zL}_{n}(\un{x}_{0}) \quad .
$$
Since the Jacobians (\ref{J4}) depend only on the $x$-coordinate and are piecewise constant, the expressions (\ref{forward})
and (\ref{backward}) take the simple form:
\bea
\overline{\zL}_{n}(i_{0})&=&- \frac{1}{n}\sum_{k=0}^{n-1}\ln J_M (i_{k}) \label{simpleforw}\\
\overline{\zL}_{n}(Q i_{n-1})&=&- \frac{1}{n}\sum_{k=0}^{n-1}\ln J_M (Q i_{(n-1)-k})= - \overline{\zL}_{n}(i_{0}) \label {simplerev}
\eea
with $i_{k}$ the region containing the point $M^k \un{x}_0$, out of the four regions
$\{ A, B, C, D \}$ and $Q=GM$, where
\be
Q A = D ~, \quad Q D = A ~, \quad Q B = B ~, \quad Q C = C \quad . \label{GM4}
\ee
Thus, the computation of $\overline{\zL}_{n}$ for the forward (respectively, time reversed) path, can be conveniently performed by keeping track only of the \textit{coarse-grained} sequences of visited regions:
\bea
\{i_k\}&=&(i_0,i_1,...,i_{n-1}) \label{path1}\\
\{Q i_{(n-1)-k}\}&=&(Q i_{n-1},Q i_{n-2},...,Q i_{0}) \label{path2} \quad  ,
\eea
rather than relying on the more detailed knowledge of the sequence of points $\{M^k\un{x}\}$  and $\{G M^{n-k}\un{x}\}$ in the phase space.\footnote{The last regions in each sequence, $i_n$ and $Q i_{-1}$, need not be taken into account, as they are unessential in the evaluation of the $\overline{\zL}_{n}$, see \cite{CKDR} for details.} A considerable amount of information, regarding the microscopic trajectory in the phase space,
is lost by passing from the phase space deterministic dynamics to the effectively stochastic process arising
from the projection of the dynamics onto the $x$-axis. Such a process is described by
a Markov jump process, which yields sequences such as those of Eqs.(\ref{path1}) and (\ref{path2}).
Nevertheless, this loss of information is irrelevant to compute the phase space contraction rate of sets of phase space trajectories.
In particular, we may disregard variations of the dynamics internal to the
single regions, as long as the resulting internal, or ``hidden'', dynamics preserve phase space volumes
and do not affect Eqs.(\ref{simpleforw}) and (\ref{simplerev}).
These observations are relevant for the stochastic descriptions of physical phenomena, which are thought
to be based on reduced (projected) dynamics of phase space deterministic dynamics. Indeed, the projected
Perron-Frobenius Equation (\ref{PF}) is a time-integrated Master equation for the probability densities $\rho_{l}$ and $\rho_{r}$ over the coarser, projected, state space defined by the chosen Markov partition.
These considerations are reminiscent of the fact that thermodynamics, to a large extent,
does not depend on the details of the microscopic dynamics, hence is consistent with many different
phase space evolutions. This is consequence of the fact that thermodynamics describes the object of
interest by a few observable quantities, i.e.\ in a space of reduced dimensionality, which can be
seen as a projection of the whole phase space.
Thus, in order to compute the quantities of interest, one may conveniently choose the detailed
microscopic dynamics which most easily represent the phenomenon under investigation.
For instance, in our idealized setting, the equilibrium dynamics and the existence of some
symmetries relating forward and reversed paths in a subspace of the phase space, may be
investigated by means of the map $M$ with $q=0$, which is reversible via the involution $G$
depicted in Fig.\ref{involM4}, and by means of the projections of $M$ on the relevant
directions.
The class of dynamics which are equivalent from a given, restricted or projected, standpoint
could include maps which are not even time reversal invariant. Indeed, one may consider microscopic
dynamics obtained from map (\ref{Map}) introducing an irreversible transformation $N$ which does
not contract nor expand phase space volumes. This can be simply done by e.g.\ letting $N$ flip the
$y$-coordinates of the phase space points of a vertical strip of width $\epsilon$ in the region $B$:
\be
\left(\begin{array}{c}
  x_{n+1} \\
  y_{n+1}
\end{array}\right)
=N
\left(\begin{array}{c}
  x_{n} \\
  y_{n}
  \end{array}\right)
  =\left\{
  \begin{array}{ll}
               \left(\begin{array}{c}
                   x_{n} \\
                   1- y_{n}
                 \end{array}\right), & \hbox{for $x \in [\tilde{x},\tilde{x}+\epsilon]$ and $y \in [0,\frac{1}{2})$} \\ \\
    \left(\begin{array}{c}
                   x_{n} \\
                   y_{n}
                 \end{array}\right) & \hbox{for $x \in [\tilde{x},\tilde{x}+\epsilon]$ and $y \in [\frac{1}{2},1]$}
  \end{array}
\right.\label{N}
\ee
cf. Fig.\ref{correction} for a graphical representation.

\begin{figure}
  \centering
   \includegraphics[width=9.5cm]{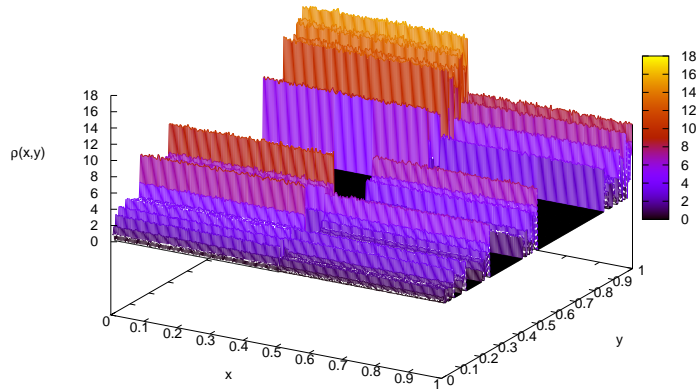}
   \caption{Result of a numerical simulation for the invariant density of the equilibrium irreversible map derived from Eqs. (\ref{Map}) and (\ref{N}) by setting $q=0$, $\ell=0.15$, $\tilde{x}=\ell$ and $\epsilon=\frac{1}{2}-\ell$, obtained by evolving a set of $2\cdot10^7$ initial conditions randomly (and uniformly) chosen on the square $[0,1]\times [0,1]$. The density $\rho(x,y)$ is not smooth along the $y$-coordinate, which is a signature of the strongly irreversible dynamics given by the map $N$ in Eq. (\ref{N}).}\label{HistB4irr}
\end{figure}

As pointed out in \cite{CKDR}, the composed map $K=N M$ is irreversible because it does not admit
an inverse. Moreover, the irreversible mechanism of the dynamics described in Eq. (\ref{N}) gives rise to an invariant measure which is fractal along the vertical direction, cf.\ Fig.\ref{HistB4irr}, and which is
strongly at variance with its reversible equilibrium counterpart shown in Fig.\ref{HistB4rev}.
Nevertheless, it is clearly seen that Eqs.(\ref{simpleforw}) and (\ref{simplerev}) still
hold true for the map $K$, despite the irreversible feature of the equations of motion. In fact, although the irreversible dynamics $K$ no longer admits an involution (hence, Eqs.(\ref{forward}) and (\ref{backward}) can
no longer be fulfilled), Eqs. (\ref{simpleforw}) and (\ref{simplerev}) remain unaltered, because $N$ maps a point
$\un{x}\in i$ into a point $\un{x}'\in i$, and it neither contracts nor expand phase space areas.
Thus, if we replace the dynamics $M$ with the dynamics $K$ and, accordingly, we take $Q=G K$, the sequences
(\ref{path1}) and (\ref{path2}) remain unaltered, since they are invariant under the action of an irreversible perturbation of the $y$-coordinate.

Now, let $\omega(i_0,n-1)\subset i_0$ and $\omega(Q i_{n-1},n-1)\subset Q i_{n-1}$ denote, respectively,
the sets of points corresponding to the forward (\ref{path1}) and to the time reversed (\ref{path2})
sequences. Trivially, the sets of points corresponding to these symbolic sequences have invariant measures
$\mu(\omega(i_0,n))$ and $\mu(\omega(Q i_{n-1},n))$, as in the reversible case. Therefore, in spite of
the irreversible modification $N$, we may say that the dynamics enjoy a form of reversibility which is
weaker than the standard reversibility in phase space, but which cannot be distinguished from that if
observed from the stochastic (reduced) viewpoint of the projections on the horizontal direction.

As a matter of fact, time reversibility is contemplated in stochastic dynamics and amounts to the requirement
that a sequence of events have positive probability if its reverse does \cite{Harris}. Because
our projected dynamics are not affected by the action of $N$, on the level of the stochastic
description, the phase space reversible dynamics of $M$ are equally {\em stochastically reversible}
as the phase space irreversible dynamics of $K$. Furthermore, the fact that $N$ reshuffles phase space
points within vertical strips of the square implies that the redistributions of mass due to the phase
space contraction rate of $M$ and to the rearrangement of phase space volumes produced by $N$ are indistinguishable on the projected horizontal space and can be quantified by the same observable $\zL$.
In particular, Eq.(\ref{SRB}) holds for both $M$ and $K$, hence the statistics of all projected observables
is the same. Only in the case that one is interested on observables which explicitly concern the $y$
coordinate would the two dynamics be distinguishable, but as long as one focuses on quantities which
do not depend on $y$, or which result from a projection on the $x$ axis, $M$ and $K$ lead to the same
conclusions: the corresponding reduced stochastic evolution is exactly the same, as in the case of
{\em felt dynamics} introduced in \cite{MR96}.

Therefore, we may state that the map $K$ with $q=0$ represents a kind of {\em equilibrium} but
{\em irreversible} dynamics. Although this might appear contradictory, it is simply explained
by the observation that the equilibrium behavior concerns the level of the horizontal projection,
which is stochastically reversible, while the irreversibility concerns the phase space.
This situation differs from that of \cite{RQ}, in which time reversible nonequilibrium systems
are recognized to be common --see, e.g.\ the standard models of nonequilibrium molecular dynamics--
while equilibrium irreversible systems are thought to be rare, as far as phase space is concerned.

Our study concerns, instead, the bridge between deterministic and stochastic-like descriptions.
In particular, we are going to show that the $x$-projection of $K$, being an {\em equilibrium} model,
satisfies the principle of detailed balance (DB), from the point of view of the stochastic dynamics,
in spite of its irreversibility in the phase space.\footnote{Note that DB is often referred to as the
{\em principle of microscopic reversibility} \cite{Tolman}, which is then to be understood as a notion
of reversibility in the reduced space, not in the full phase space.}

This naturally connects with the distinction between {\em relevant} and {\em irrelevant} coordinates
which underlies the statistical mechanics reduction of deterministic descriptions in the phase space to
stochastic descriptions, which typically concern the one-particle space \cite{Onsager,RV}. Clearly, our dynamical system is too simple to allow a physically meaningful distinction between relevant and irrelevant variables. Therefore, we merely focus on one of them (the $x$-coordinate) and regard the other as irrelevant  (the $y$-coordinate), without implying that one subspace of our subspace is endowed with any special meaning.

\section{From phase space to detailed balance}
\label{sec:sec2}

In this Section we introduce the notion of \textit{detailed balance in the phase space} (PSDB) which constitutes a strong concept of equilibrium dynamics. We will show that the standard notion of DB descends
from PSDB through a projection of the phase space dynamics onto a suitable subspace. It will also be shown that the DB condition is insensitive to the reversibility properties of the full phase space equilibrium dynamics.
Let
$$
 \un{x}_{n+1}=M \un{x}_{n}~, \quad \un{x} \in \mathcal{U} \label{propag}
$$
be the microscopic dynamics, where $M$ is time reversal invariant, with involution $G$. For sake of
simplicity, we deal with discrete time, $t \in \mathbb{Z}$, but flows $S^t$, $t \in \mathbb{R}$
could be treated similarly.

\begin{figure}
\begin{center}
  \includegraphics[width=9cm, width=0.70\textwidth]{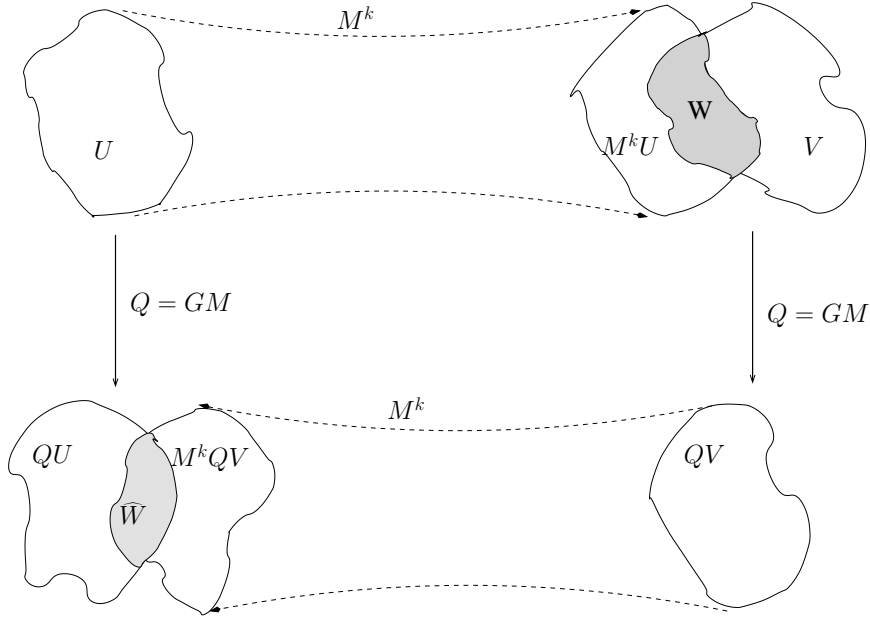}\\
  \caption{Set of points, in phase phase, belonging to the sets $U$ and $V$ in the forward (upper path) and in the time reversed (lower path) trajectories}\label{GDB}
   \end{center}
\end{figure}

Consider two sets in phase space, $U,V \subset\mathcal{U}$. Let
$W=M^{k}U\cap V = \{\un{x} \in V : M^{-k} \un{x} \in  U \}$ be the set of final points of
trajectory segments starting in $U$, which fall in $V$ after $k$ iterations of $M$,
and let $M^{-k}W$ be the corresponding set of initial conditions. Take an invariant
measure $\mu$ for $M$, so that $\mu(W) = \mu(M^{-k}W)$, and let $\zL_k(\un{x})$ be
the phase space contraction along the trajectory starting in $\un{x} \in M^{-k}W$.
Its time reverse trajectory, which gives rise to the opposite phase space contraction,
cf.\ Eqs.(\ref{backward}) and (\ref{simplerev}), starts in
$\hat{\un{x}} = G M M^k \un{x} \in G M W$. Call
$$
\widehat{W}= M^{k} G M W = M^{k} G M \left( M^{k}U \cap V \right) =
G M U  \cap M^{k} G M V
$$
the set of final conditions of all trajectories which are time reverses of those ending in
$W$, cf. Fig. \ref{GDB}, where the second equality comes from the definition of $W$ and the third equality follows from the property $M^{k} G= G M^{-k}$ of the involution $G$. The measure of this set is given by:
\bea
\mu(\widehat{W})&=&\mu(M^{k}G M V\cap G M U) \nonumber\\
&=& \mu(M^{k} M^{-1} G V\cap M^{-1} G U)=\mu(M^{k-1}( G V\cap M^{-k} G U)) \nonumber\\
&=&\mu(M^{k-1}( G V\cap G M^{k} U))= \mu(G V\cap G M^{k} U)\nonumber\\
&=& \mu(G (V\cap M^{k} U))=\mu(GW) \label{GW}
\eea

We define the {\em phase space detailed balance} (PSDB) as the condition for which the probability of having
opposite phase space contractions are equal:
\be
\mu(M^{-k} W) = \mu(M^{-k} \widehat{W}) \quad \mbox{ i.e. } \quad
\mu(W)=\mu(\widehat{W})  \label{PSDB-0}
\ee
Becasue of Eq.(\ref{GW}), this condition may also be written as
\be
\mu(W)=\mu(GW)  \label{PSDB}
\ee
which is to say that PSDB requires the $M$-invariant measure $\mu$ to be also $G$-invariant,
since $W$ may be any subset of $\mathcal{U}$. Calling a function $\Phi$ odd with respect to time
time reversal if $\Phi (G\un{x})=-\Phi (\un{x})$ for all $\un{x} \in \mathcal{U}$,
and defining equilibrium the situation in which the mean value of all such odd observables vanishes,
we obtain that Eq.(\ref{PSDB}) implies equilibrium. Indeed, take any $\Phi$ which is odd with respect
to the time reversal and choose $W = \{ \un{x} \in \mathcal{U} : \Phi(\un{x}) > 0 \}$, so that
$G W = \{ \un{x} \in \mathcal{U} : \Phi(\un{x}) < 0 \}$ and $\Phi(\un{x}) = 0$ for all
$\un{x} \in \mathcal{U} \setminus (W \cup GW)$. Then the following holds:
\bea
\langle \Phi \rangle &=&\int_W \Phi(\un{x}) d\mu(\un{x})+\int_{GW}\Phi(\un{x}) d\mu(\un{x})
\nonumber\\
&=& \int_W \Phi(\un{x}) d\mu(\un{x})+\int_{W}\Phi(\un{G y}) J_G(\un{y}) d\mu(\un{y})
\\
&=& \int_W \Phi(\un{x}) d\mu(\un{x})-\int_{W}\Phi(\un{y})  d\mu(\un{y})= 0 \nonumber
\label{nocurr}
\eea
where $J_G=1$. In particular, the PSDB implies that the average of the phase space
contraction rate $\Lambda(\un{x})=\log J_F(\un{x})^{-1}$ vanishes, as required for equilibrium
in Sec.\ref{sec:sec2}.
\begin{figure}
\begin{center}
  \includegraphics[width=9cm]{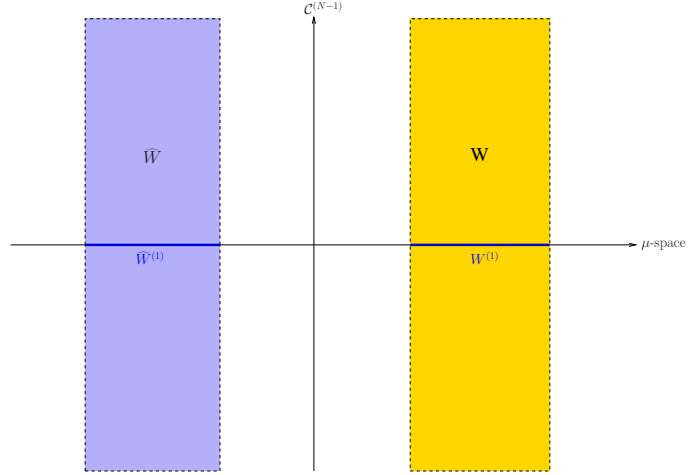}\\
  \caption{The sets $W^{(1)}$ and $GW^{(1)}$, after reduction from phase space to the $\mu$-space}\label{proj}
   \end{center}
\end{figure}

Is there any relation between PSDB and the standard DB, which implies equilibrium on the
level of the projected dynamics (the caricature of the one particle or $\mu$-space)?\footnote{Observe that
some authors distinguish DB dynamics from detailed balance steady state, which proves to be a
convenient tool in the analysis of stochastic processes \cite{CMW}. In this case, a given evolution
law is called DB dynamics if its steady state is a DB state.}
To derive standard DB, let us eliminate the irrelevant coordinates, which only contribute to noise,
by projecting Eq.(\ref{PSDB}) on the subspace of {\em relevant} coordinates, which we call $\mu$-space.
This can be done for sets of the form
$W=W^{(1)} \times \mathcal{C}_{n-1}$ and $\widehat{W}=\widehat{W}^{(1)}\times \widehat{\mathcal{C}}_{n-1}$,
where $W^{(1)}$ and $\widehat{W}^{(1)}$ denote the projections of the sets $W$ and $\widehat{W}$ onto the $\mu$-space and $\mathcal{C}_{n-1},\widetilde{\mathcal{C}}_{n-1}$ span the remaining, noisy, space, cf.\ the hyper-cylinders illustrated in Fig. \ref{proj}, where $W^{(1)}$ and $\widehat{W}^{(1)}$ denote, respectively, the event on the $\mu$-space.
Then, if we denote by $\mu^{(1)}$ the measure in the $\mu$-space induced from the invariant measure
on the phase space, we have that:
\be
\mu^{(1)}\left(W^{(1)}\right)=\int_{\mathcal{C}_{n-1}}\mu(d\un{x}) ~, \quad
\mu^{(1)}\left(\widehat{W}^{(1)}\right)=\int_{\widehat{\mathcal{C}}_{n-1}}\mu(d\un{x}) \label{mu1G}
\ee
Detailed balance holds if
\be
\mu^{(1)}\left(W^{(1)}\right)=\mu^{(1)}\left(\widehat{W}^{(1)}\right)~, \quad \mbox{or } \quad
\mu^{(1)}\left(W^{(1)}\right)=\mu^{(1)}\left(GW^{(1)}\right) \label{DB}
\ee
Thus, PSDB implies the standard DB, since DB amounts to the condition of PSBD restricted to
sets $W$ and $\widehat{W}$ of the form here introduced. Moreover, the projection procedure
typically smoothes out singularities, hence the induced invariant measure $\mu^{(1)}$ usually
is regular and has an invariant density $\rho^{(1)}$.\footnote{For instance, the map (\ref{Map})
with $h=\frac{1}{2}-2\ell$ which, for arbitrary $\ell$, is dissipative but still equipped with a
projected invariant density which is smooth along the unstable manifold, \cite{CKDR}.}

In the derivation of the stochastic description as a projection of some deterministic phase
space dynamics, the DB condition (\ref{DB}) is usually assumed to be the consequence of the
time reversal invariance of the microscopic dynamics, once equilibrium is reached. In our
investigation this amounts to require the existence of the involution $G$, defined in the
phase space. We are now going to see that this requirement may be relaxed in simple cases,
such as those of the irreversible dynamics $K$ discussed above. Using the same notation,
Eq.(\ref{PSDB}) can then be written as:
\be
\mu(Mi\cap j)=\mu(M G M j \cap G M i) \label{simple2}
\ee
where $\mu(Mi\cap j)$ is the conditional probability of being in the region $j$ one time step
after having been in region $i$, with $i,j \in \{A,B,C,D\}$, our finite state space.
The quantity $\mu(Mi\cap j)$ may be rewritten as $\mu(M i \cap j) = p(j|i)=p_{ij}$, with notation
reminiscent of stochastic descriptions. The $p_{ij}$'s then constitute the elements of the transition
matrix:
\be
 P=\left(
     \begin{array}{cccc}
       0 & 0 & \frac{1}{2} & \frac{1}{2} \\
       2\ell & 1-2\ell & 0 & 0 \\
       0 & 0 & \frac{1}{2} & \frac{1}{2} \\
       2\ell & 1-2\ell & 0 & 0 \\
     \end{array}
   \right) \label{P}
 \ee
which defines a {\em stochastic} process for the dynamics of the state in $\{A,B,C,D\}$.
Application of the Perron-Frobenius Theorem to the matrix $P$ reveals the existence of one left eigenvector of $P$, associated with the eigenvalue $\lambda=1$ (i.e. the geometric multiplicity of such eigenvalue is $1$), which, thus, implies the existence of a unique (coarse-grained) steady state measure $(\mu_A , \mu_B , \mu_C , \mu_D)$, where:
\be
\mu_i=
\left\{
  \begin{array}{ll}
    \dfrac{2\ell}{1+4\ell}, & \hbox{if ~ $i=A,C,D$;} \\ \\
    \dfrac{1-2\ell}{1+4\ell}, & \hbox{if ~ $i=B$.}
  \end{array}
\right.\label{mu}\quad .
\ee
Eq. (\ref{mu}) highlights the fact that the measure in the full phase space is uniform along the stable manifold and piecewise constant along the unstable one, as previously illustrated in Fig.\ref{HistB4rev}. This comes from the fact that, since the regions $\{A,B,C,D\}$ are determined only by the $x$-coordinate, the stochastic transition matrix $P$ does not depend on the dynamics along the stable manifold, and, hence, $p_{ij}=(J_j^u)^{-1}$. As a result, the invariant measure in Eq. (\ref{mu}) depends only on $x$ and is constant on $y$.

\begin{figure}
\begin{center}
  \includegraphics[width=8.5cm]{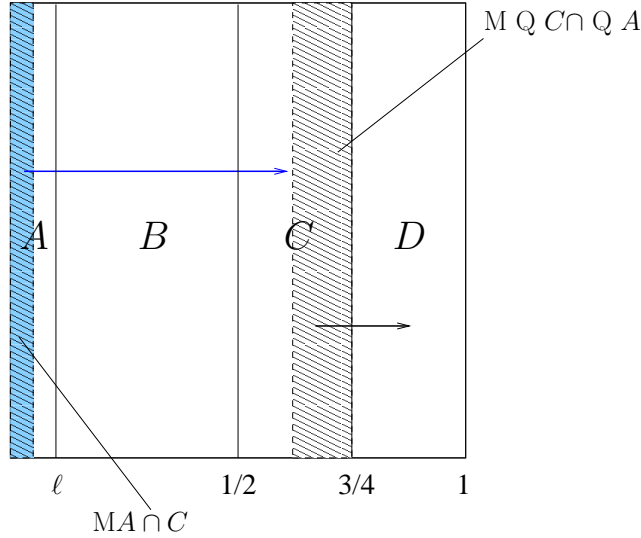}\\
  \caption{Sets undergoing the forward path ($A\rightarrow C$, blue arrow) and the time reversed path ($QC \rightarrow QA$ black arrow) for the map (\ref{Map}) with $q=0$.}\label{jump}
   \end{center}
\end{figure}
Consider, for instance, the one-step transition $A \rightarrow C$, whose probability is
the measure of $W=M A \cap C$, which equals that of the set of the corresponding initial
conditions $M^{-1} W = A \cap M^{-1} C$. The probability of the reverse transition
$Q C = C \rightarrow Q A = D$, where we have recalled the relations (\ref{GM4}), is
the measure of the set $\widehat{W}=M Q C \cap Q A$, cf.\ Fig.\ref{jump}.
Since the sets $W$ and $\widehat{W}$ span the whole range $[0,1]$ in the vertical direction,
as in Fig. \ref{proj}, the measures $\mu(W)$ and $\mu(\widehat{W})$ can be calculated just in
terms of the transition probabilities (\ref{P}) and of the projected invariant measures
(\ref{mu}). The result is
\be
\mu(MA\cap C)=\mu_A p_{AC}=\frac{\ell}{1+4\ell}=\mu_C p_{CD}=\mu(M C\cap D) \quad \
\ee
for any $\ell$. Hence, DB holds as expected, because we are dealing with an equilibrium case,
although the underlying dynamics is irreversible.
It is now interesting check what happens when the microscopic dynamics is pulled out of equilibrium.
For instance, consider $K=N M$, with $N$ as in (\ref{N}) and $M$ the time reversible mapping
(\ref{Map}), with $q=\frac{1}{2}-2\ell$.
For this map $M$, we may consider an involution $G$ which is consistent with the following
equalities (cf.\ Eqs.(36) in \cite{CKDR}):
\be
Q A = A ~, \quad Q B = C ~, \quad Q C = B ~, \quad Q D = D \label{GM3}
\ee
hence which differs from the $G$ in Eq.(\ref{G}). Then, the time reverse of the transition
$A \rightarrow C$ is given by $Q C \rightarrow Q A$ i.e.\ $B \rightarrow A$, and we get:
\bea
\mu(MA\cap C)&=&\mu_A p_{AC}=\frac{\ell}{1+4\ell} \label{old1}\\
\mu(MB\cap A)&=&\mu_B p_{BA}=\frac{(1-2\ell)^2}{1+4\ell} \label{old2}
\eea
which shows that (\ref{old1}) and (\ref{old2}) do not coincide, and that PSDB and DB are violated for
$q\neq0$, i.e.\ outside the equilibrium defined via the chosen $G$.
In this case, only $\ell=\frac{1}{4}$ leads to the equilibrium state which is,
in addition, microcanonical\footnote{The map featuring $q=\frac{1}{2}-2\ell$ attains equilibrium for $q=0$, which gives $\ell=\frac{1}{4}$, corresponding, as discussed in Ref.\cite{CKDR}, to a microcanonical equilibrium distribution.} .

\section{The fluctuation relation and nonequilibrium response}
\label{sec:sec3}

The fluctuation relation for $\Lambda$, the $\Lambda$-FR, originally proposed by
Evans, Cohen and Morriss \cite{ECM}, and developed by Gallavotti and Cohen \cite{Gall}, concerns
the statistics of the mean phase space contraction rate $\overline{\zL}_{n}$, over the steady state
ensemble of phase space trajectory segments of a large number of steps, $n$.
Equivalently, it concerns the statistics of $\overline{\zL}_{n}$, computed over segments of a unique
steady state phase space trajectory, broken in segments $\{{x}_1,...,{x}_n\}$ of length $n$.

The dynamics are called dissipative if
$$
\langle \zL \rangle = \int_\mathcal{U} \zL(\un{x}) ~ \mu({\rm d} \un{x}) > 0
$$
where $\langle \zL \rangle$ is the steady state mean of $\zL$, i.e.\
it is computed with respect to the natural measure $\mu$ on $\mathcal{U}$.
It is convenienet to introduce the dimensionless phase space contraction rate $e_{n}={\overline{\zL}_{n}}/{\langle \zL \rangle}$ because its range,
$[-\zL_{\mbox{max}}/\langle \zL \rangle,\zL_{\mbox{max}}/\langle \zL \rangle]$,
does not change with $n$, while the values taken by $e_n$ tend to more and more densely
fill it as $n$ grows. Further, we denote by $\pi_n(B_{p,\zd})$
the probability that $e_n$, computed over a segment of $n$-steps of a typical trajectory, falls
in the interval $B_{p,\zd}=(p-\zd,p+\zd)$, for some fixed $\zd>0$. In other words, one may
write
\be
\pi_n(B_{p,\zd}) = \mu(\omega_{\Lambda,n}) ~, \quad \mbox{where } \quad
\omega_{\Lambda,n}=\{\un{x}\in \mathcal{U}: e_{n}\in B_{p,\zd}\}
\ee
For growing $n$, $\pi_n$ peaks around the mean value $\langle e_n \rangle =1$, but fluctuations about
this mean may occur with positive probability at any finite $n$. In particular, under certain
conditions, \cite{Gall,Gawedzky,Ellis,SearlRonEvans}, $\pi_n$ obeys a large deviation principle
with a given rate functional $\zeta$, in the sense that the limit
\be
\lim_{n\rightarrow \infty} \pi_n(B_{p,\zd})= e^{-n[\zeta(p) + \epsilon_\delta]}
\label{largedev}
\ee
exists, with $\epsilon_\delta \le \delta$. In particular, if the support of the
invariant measure is the whole phase space $\mathcal{U}$, time reversibility
guarantees that the support of $\pi_n$ is symmetric around $0$, and one can
consider the ratio
$$
\frac{\pi_n(B_{p,\zd})}{\pi_n(B_{-p,\zd})} ~ .
$$
In our case, this ratio equals the ratio of the measures of a pair of sets conjugated by time
reversal, as in Eqs. (\ref{simpleforw}) and (\ref{simplerev}). Then, the validity of the
$\zL$-FR means that there exists $p^* > 0$ such that
\begin{equation}
p - \zd \le \lim_{n \to \infty} \frac{1}{n \langle \zL \rangle} \log
\frac{\mu(\{ x : e_n(x) \in B_{p,\delta}\})}
{\mu(\{ x : e_n(x) \in B_{-p,\delta} \})} \le p+\zd
\label{prethm}
\end{equation}
if $|p| < p^*$ and $\zd > 0$.

So far, the proofs of this and other FR's appeared in the literature, notably those for the
fluctuations of the Dissipation Function $\Omega$, \cite{RonMejia,SearlRonEvans,MoRo}, rely
on the existence of an involution representing time reversal in phase space, while they rely
on the principle of microscopic reversibility in the state space of stochastic processes. So,
whatever the context, the relevant notion of time reversibility has always been used.
Therefore, if time reversibility is broken, but is broken as in the case of the map $K$,
which enjoys a weaker form of reversibility requiring only the existence of the pairs of
conjugate trajectories (\ref{path1}) and (\ref{path2}), the $\Lambda$-FR
should remain valid. Indeed, the existence of this weaker reversibility is consistent with
the principle of microscopic reversibility, adopted in the stochastic approach by e.g.\
Lebowitz and Spohn, \cite{LeboSpohn}.

Let us then investigate the validity of the $\Lambda$-FR for the deterministic model $K=N M$
with $N$ given by Eq.(\ref{correction}) and $M$ by Eq.(\ref{Map}), which may or may not lead to
equilibrium, depending on the value of the parameter $q$. As discussed in Sec. \ref{sec:sec1},
$K$ is irreversible, $M$ is reversible and the mapping $Q=GK$, appearing in the definition of
the time reversed path (\ref{path2}), is properly defined in both cases.
Take $q=\frac{1}{2}-2\ell$, with $\ell \neq \frac{1}{4}$, consistently with \cite{CKDR}. Then,
the $\Lambda$-FR may be written as:
\be
p - \zd \le \lim_{n \to \infty} \frac{1}{n \langle \zL \rangle} \log
\frac{\mu(\omega(i_0,n-1))}
{\mu(\omega(Q i_{n-1},n-1))} \le p+\zd
\label{prethm3}
\ee
with $|p| < p^*$, for some $p^* > 0$ and any $\zd > 0$.
To prove Eq.(\ref{prethm3}) for the map $K=N M$, one must compute the invariant probability measure
$\pi_n$ and let $n$ grow without bounds.

\noindent
{\em This is guaranteed by the proof of the validity
of the $\Lambda$-FR given for $M$ in Ref.\cite{CKDR}, which only relies on the invariant measure in
the projected space. }

\begin{figure}
\begin{center}
  \includegraphics[width=10.5cm]{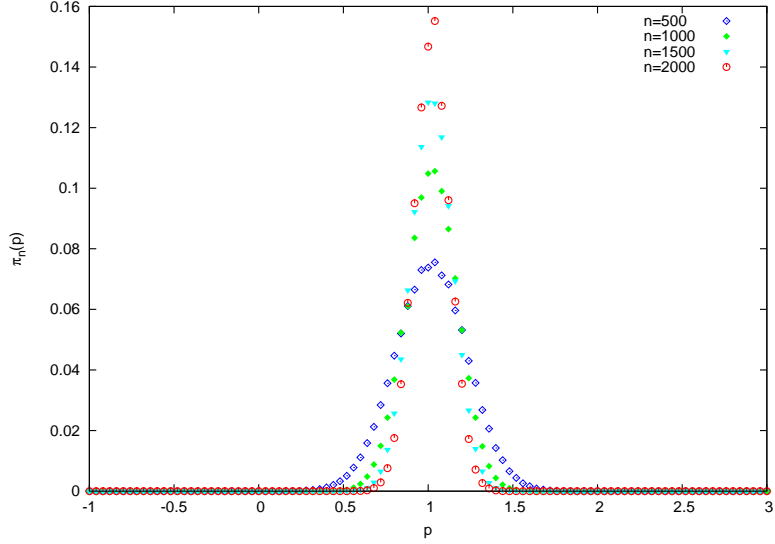}\\
  \caption{Probability measure $\pi_n(B_{p,\zd})$}\label{prob}
   \end{center}
\end{figure}

\begin{figure}
\begin{center}
  \includegraphics[width=10.5cm]{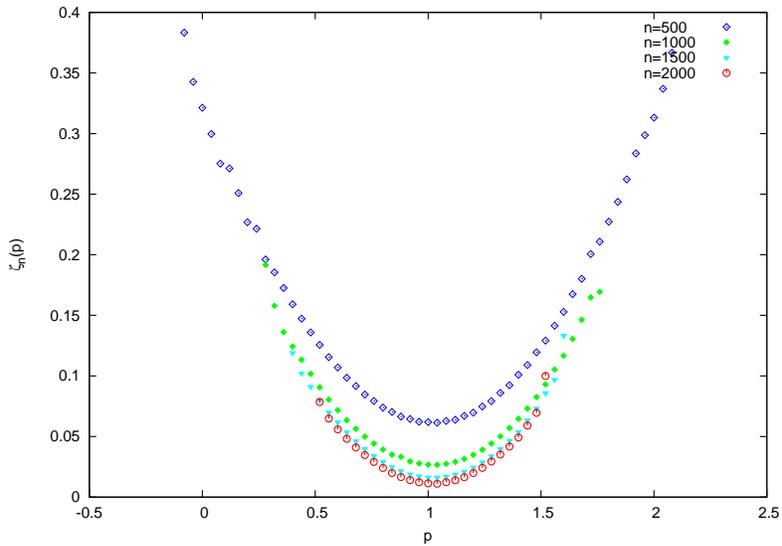}\\
  \caption{Rate functional $\zeta_n$ associated with $\pi_n$for different values of $n$. As expected from theory, the curves $\zeta_n$ are expected to move downwards for growing $n$, so that, in the
$n\rightarrow \infty$ limit, $\zeta(p)$ intersects the horizontal axis only in $p=1$.}
\label{ratefunct0}
   \end{center}
\end{figure}

\begin{figure}
\begin{center}
  \includegraphics[width=10.5cm]{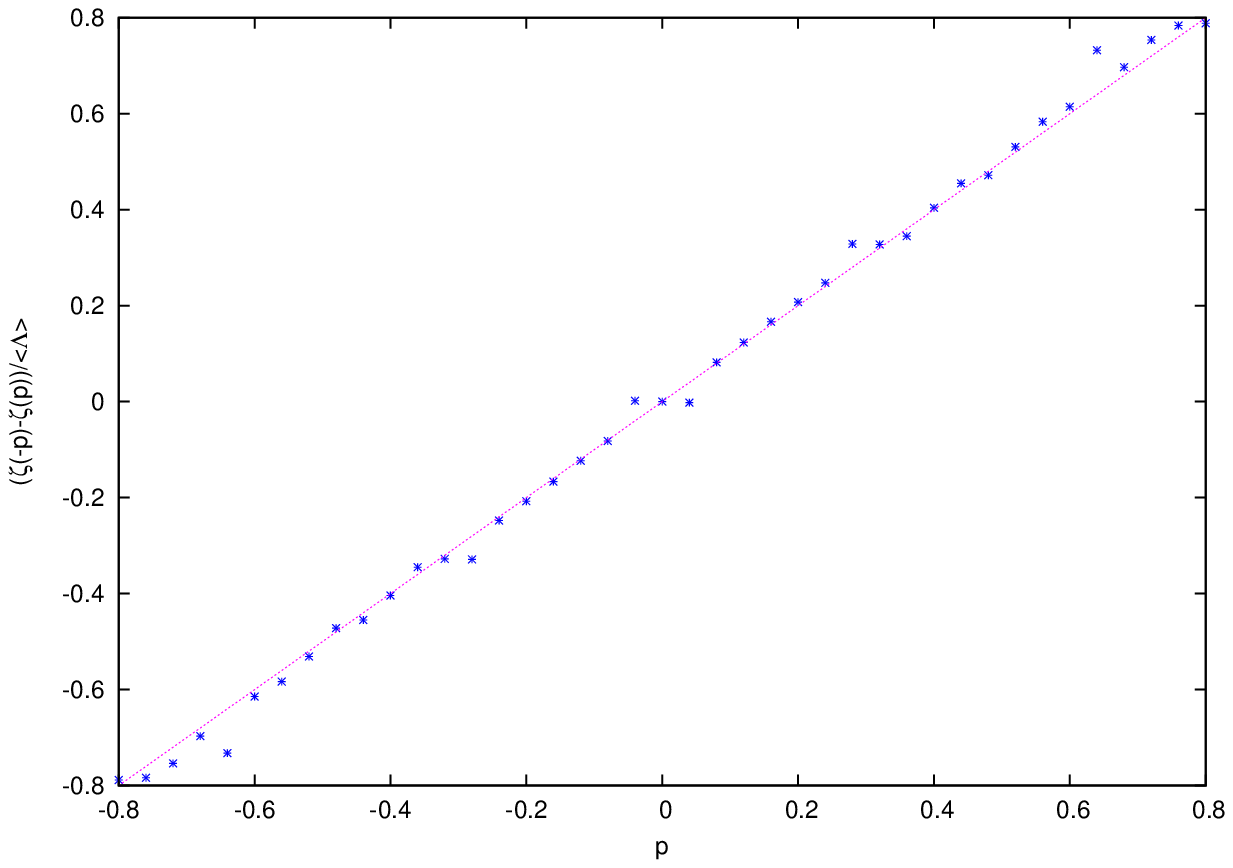}\\
  \caption{Check of the $\Lambda$-FR, Eq. (\ref{prethm3}) for the irreversible dynamics dictated by the map $K=NM$. \textit{Blue points}: results of the numerical simulations with $n=2\cdot10^2$.}\label{FR}
   \end{center}
\end{figure}

We illustrate this result by means of numerical simulations, which we
have performed for different values of $n$, with
$\ell=0.15$, cf. Fig.\ref{prob}. The numerical simulations, shown in Fig. \ref{ratefunct0}, show how the large deviation
rate functional $\zeta$ is generated, and that it is smooth and strictly convex in
the whole range of observed fluctuations, as required by the theory of the $\Lambda$-FR.
In Fig.\ref{FR} we also plotted the expression in the center of Eq.(\ref{prethm3}), which is
consistent with the validity of the $\Lambda$-FR.
Once the $\Lambda$-FR is proven to hold for the map $K$, one may be tempted to assess the validity of the
Green-Kubo formulas as well as of the Onsager reciprocal relations, by following e.g.\ the
strategies of \cite{SRE,Gall2,RoCoh,Ladek}, in the limit of small external drivings, as summarized
in \cite{Bett,RonMejia}. To this end, we consider $\Lambda$ as the entropy production
rate,\footnote{This identification must be done {\em cum grano salis}, as explained in e.g.\
\cite{SRE,RCNon2000,RTV00,CRPhysA202,RKbook}.} and we briefly summarize the the argument, for
sake of completeness.

The main steps are the following \cite{Bett,Gall2}:
\begin{itemize}
\item
assume that the system is subjected to $k$ fields $F=(F_1,F_2,..,F_k)$, that $\Lambda$ vanishes
when all drivings vanish and that
\be
\Lambda(\un{x})=\sum_{\ell=1}^{k}F_\ell J_\ell^0(\un{x})+O(F^2)
\ee
which defines the currents $J_r^0$, which are proportional to the forces $F_r$.
\item
The decay of the $\Lambda$ autocorrelation function required for the $\Lambda$-FR to hold,
leads to the following expansion for the rate function $\zeta$:
\be
\zeta(p)=\frac{\langle \zL \rangle^2}{2C_2}(p-1)^2+O((p-1)^3F^3)
\ee
where $C_2$ is related to the time autocorrelation of $\zL$. In other words, the rate functional is quadratic
for small deviations from the mean $p=1$, in accord with the Central Limit Theorem.
\item
Introduce the nonlinear currents as $J_\ell(\un{x})=\partial_{F_\ell}\Lambda(\un{x})$, and the
transport coefficients $L_{\ell r}=\partial_{F_r}\langle J_\ell \rangle |_{F=0}$. Then, one obtains:
\be
\langle \Lambda \rangle=\frac{1}{2}\sum_{\ell,r =1}^{k}(\partial_{F_r}\langle J_\ell \rangle +\partial_{F_\ell}\langle J_r \rangle)|_{F=0} F_\ell F_r=
\frac{1}{2}\sum_{\ell,r =1}^{k}(L_{\ell r}+L_{r\ell})F_\ell F_r \label{GK}
\ee
to second order in the forces.
\item
Equation (\ref{prethm3}) implies $\langle \Lambda \rangle=\frac{C_2}{2}$.
Thus, equating the latter expression with Eq.(\ref{GK}) and by considering $(L_{\ell r}+L_{r\ell})/2$ with $\ell=r$, one obtains the Green-Kubo relations.
\end{itemize}

\begin{figure}
\begin{center}
  \includegraphics[width=10.5cm]{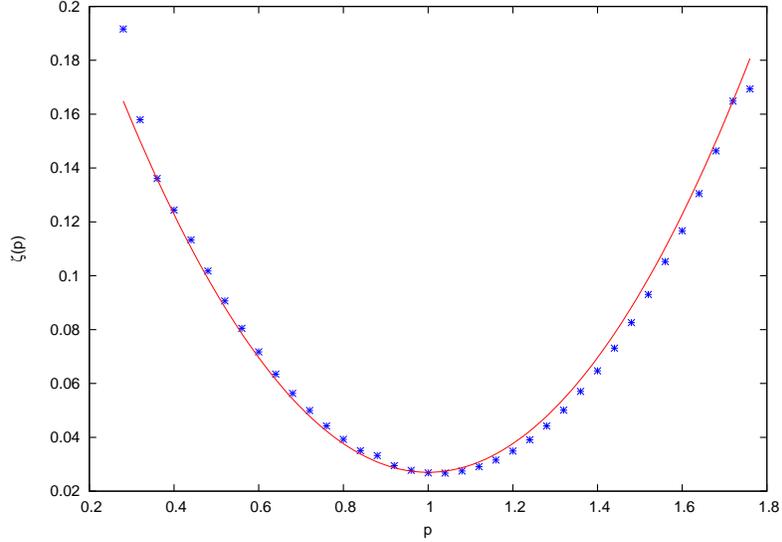}\\
  \caption{Rate functional $\zeta(p)$ associated to $\pi_n(B_{p,\zd})$. \textit{Blue points}: results of the numerical simulation with $n=10^3$. \textit{Red line}: Fitting of numerical data with the parabola $a(x-1)^2+b$, with parameters $a=2.66\cdot10^{-1}$ and $b=2.68\cdot10^{-2}$.}\label{ratefunct}
   \end{center}
\end{figure}

In our case, the rate functional is  clearly quadratic, as shown by our simulations of the dynamics of $K$.
In particular, the red quadratic curve in Fig.\ref{ratefunct} reproduces nicely the behavior of the numerical
data for the rate functional $\zeta_n$, corresponding to trajectory segments of $n=200$ steps.
The necessity for a parameter $b \neq 0$ in the parabola is due the finiteness of $n$: indeed
$b \rightarrow 0$ when $n \rightarrow \infty$.
However, in spite of the validity of the $\Lambda$-FR for the irreversible map $K$,
which entails that the irreversible map behaves to some extent equivalently to the reversible map
$M$, the argument of \cite{Gall98} leading to the Green-Kubo relations cannot be reproduced here. In fact,
it relies on the differentiability of the SRB measure as well as on the reversibility of the microscopic
dynamics, which are both violated in the case of $K$.
Alternatively, one may think of deriving linear response from the $\Lambda$-FR through the approach
of Ref.\cite{SRE}  (SRE, hereafter), which does not explicitly require the differentiability of the
invariant measure and the reversibility of the dynamics. In particular, SRE deals with a Nos\'e-Hoover
thermostatted $N$-particle system and obtains:
\be
n \sigma^2_{\overline{\zL}_{n}}(F_e)=\left(\frac{F_e V}{2 K_0}\right)^2
\left[\frac{2L(F_e)k_B T}{V}+O\left(\frac{F_e^2}{n N}\right)\right] \label{R1}
\ee
for the variance $\sigma^2_{\overline{\zL}_{n}}$ of $\overline{\zL}_{n}$, provided $\zL$ can be identified
with the dissipation function $\Omega$. Here, where $K_0$ is the target kinetic energy for the thermostatted
particles, which corresponds to
the inverse temperature $\beta$, $N$ is the number of particles, $V$ is the volume
and $F_e$ denotes the external force (i.e. the bias) acting on the system. The quantity $L(F_ e)$ is defined by
\be
L(F_e)=\beta V\int_0^\infty dt \langle (\Psi(t)-\langle \Psi \rangle) (\Psi(0)-\langle \Psi \rangle)\rangle \label{R2}
\ee
and $L(0)=\lim_{F_e \rightarrow 0}L(F_e)$ is the linear transport coefficient.
The derivation of the Green-Kubo formulae is completed by comparing (\ref{R1}) with the relation
$\langle \Lambda \rangle= \langle \Psi \rangle F_e = \frac{n}{2} \sigma^2_{\overline{\zL}_{n}}(F_e)$,
which is implied by the FR, which yields:
\be
L(0)=\lim_{F_e \rightarrow 0}\frac{\langle \Psi \rangle}{F_e}=\beta V\int_0^\infty dt \langle \Psi(t)\Psi(0)\rangle \label{R3}
\ee

In our framework of simple dynamical systems, the ``current'' $\Psi$ could be defined as:
\be
\Psi(\un{x})=\left\{\begin{array}{cc}
                 0 & \quad \text{for $\un{x} \in A, D$} \\ \\
                 1 & \quad \text{for $\un{x} \in B$} \\ \\
                -1 & \quad \text{for $\un{x} \in C$}
               \end{array}\right.
 \; ,\;  \label{J}
\ee
which implies an average current $\langle \Psi \rangle = {(1-4 \ell)}/{(1+4 \ell)}$, cf.\ Eq.(41) in Ref.\cite{CKDR}, where the role of the external force $F_e$ was played by the \textit{bias} $b=2-{1}/{(1-2\ell)}$. \footnote{Adopting the definition of current in Eq. (\ref{J}) for the map of Sec. \ref{sec:sec2}, one
finds a ``positive current'' for expanding phase space volumes. This is different from the case of the most common deterministically thermostatted dynamics, typical of nonequilibrium molecular dynamics. However, there is no general principle which imposes phase space contraction in presence of positive currents. Indeed, one may consider models without any phase space variations, or models with phase space expansion for positive currents (see e.g. \cite{BR} for certain parameter values). Moreover, in our highly idealized model, the observable $\Psi$ can be defined differently, so that it may take whatever values one likes.}

\begin{figure}
   \centering
   \includegraphics[width=10.5cm]{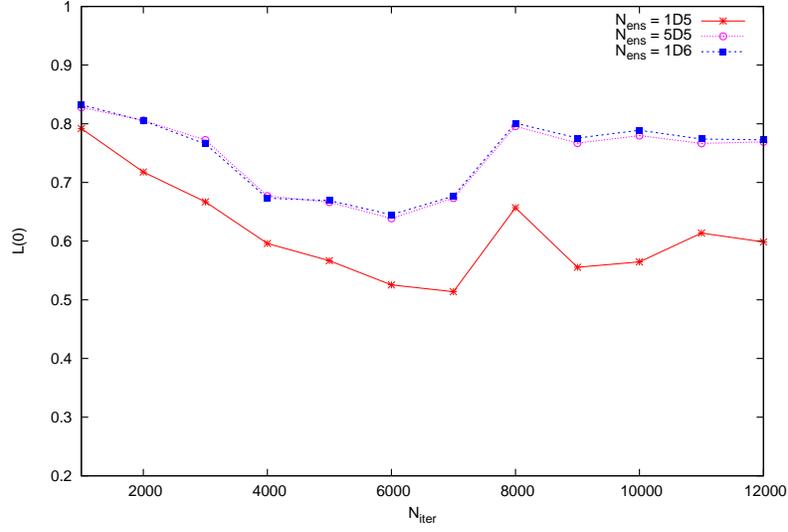}
   \caption{Numerical computation of $L(0)$ from Eq.(\ref{M1}), for the irreversible map $K=M N$ derived from Eqs. (\ref{N}) and (\ref{Map}), with $q=\frac{1}{2}-2\ell$, for different values of $N_{ens}$.}\label{conv1}
\end{figure}

Nevertheless, following SRE may be problematic, as $N$ is required to be large enough in order to derive (\ref{R3}) from (\ref{R1}), something which cannot be granted in low-dimensional systems as the map under consideration.
Indeed, our numerical simulations reveal that an interesting scenario arises in the computation of the
quantity $L(F_e)$, which may be conveniently approximated by:
\be
L(F_e)\simeq\frac{1}{N_{ens}}\sum_{k=0}^{(N_{iter}-1)}\sum_{j=1}^{N_{ens}}
\left[\Psi(\un{x}_{k}^{(j)})\Psi(\un{x}_{0}^{(j)})-\langle \Psi \rangle^2 \right] \label{M1}
\ee
where the upper limit of the integral in (\ref{R2}) is replaced by $N_{iter}$ and the correlations
are computed over an ensemble of fully decorrelated initial conditions $\{\un{x}_{0}^{(j)}\}$,
with $j=1,...,N_{ens}$, picked at random so that they occur in the ensemble with the frequency
corresponding to the natural invariant measure of the dynamical system.
Then, the coefficient $L(0)$ may be computed from Eq. (\ref{M1}) by considering the limit of vanishing bias, i.e. by taking $\langle \Psi \rangle=0$ and a microcanonical \footnote{See footnote at the end of Sec.\ref{sec:sec2}.} equilibrium ensemble of initial conditions $\{\un{x}_0^{(j)}\}$.
The values $N_{ens}$ and $N_{iter}$ must be chosen with care, in order to
guarantee the convergence of the sums in (\ref{M1}), cf. Fig.\ref{conv1}.

\begin{figure}
  \centering
   \includegraphics[width=10.5cm]{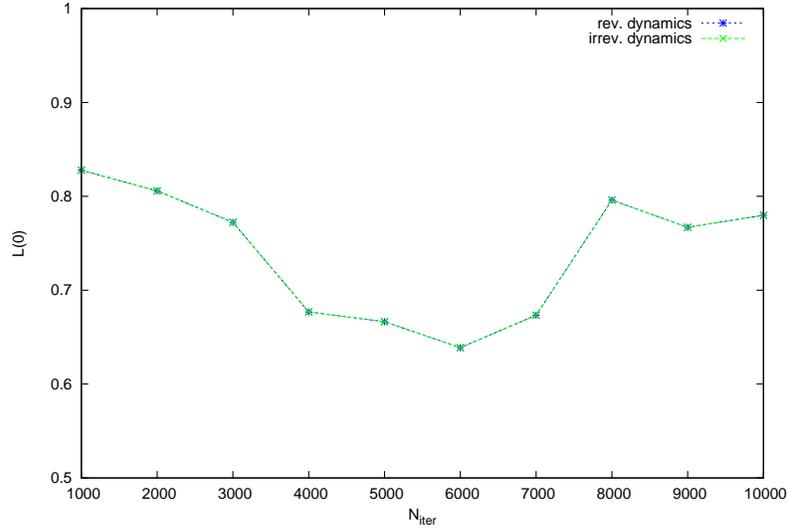}
   \caption{Computation of $L(0)$ for the irreversible, $K=M N$, and the corresponding reversible map, $K=M$, with $N_{ens}=5\cdot 10^5$. The irreversibility does not affect the convergence in Eq.(\ref{M1}).}\label{conv2}
\end{figure}

Our simulations show that the irreversible character of the dynamics, which affects the ``irrelevant"
variable $y$, has no influence on the convergence of the Green-Kubo formula Eq.(\ref{M1}), cf.\ Fig.\ref{conv2}.

\begin{figure}
\centering
\includegraphics[width=6cm]{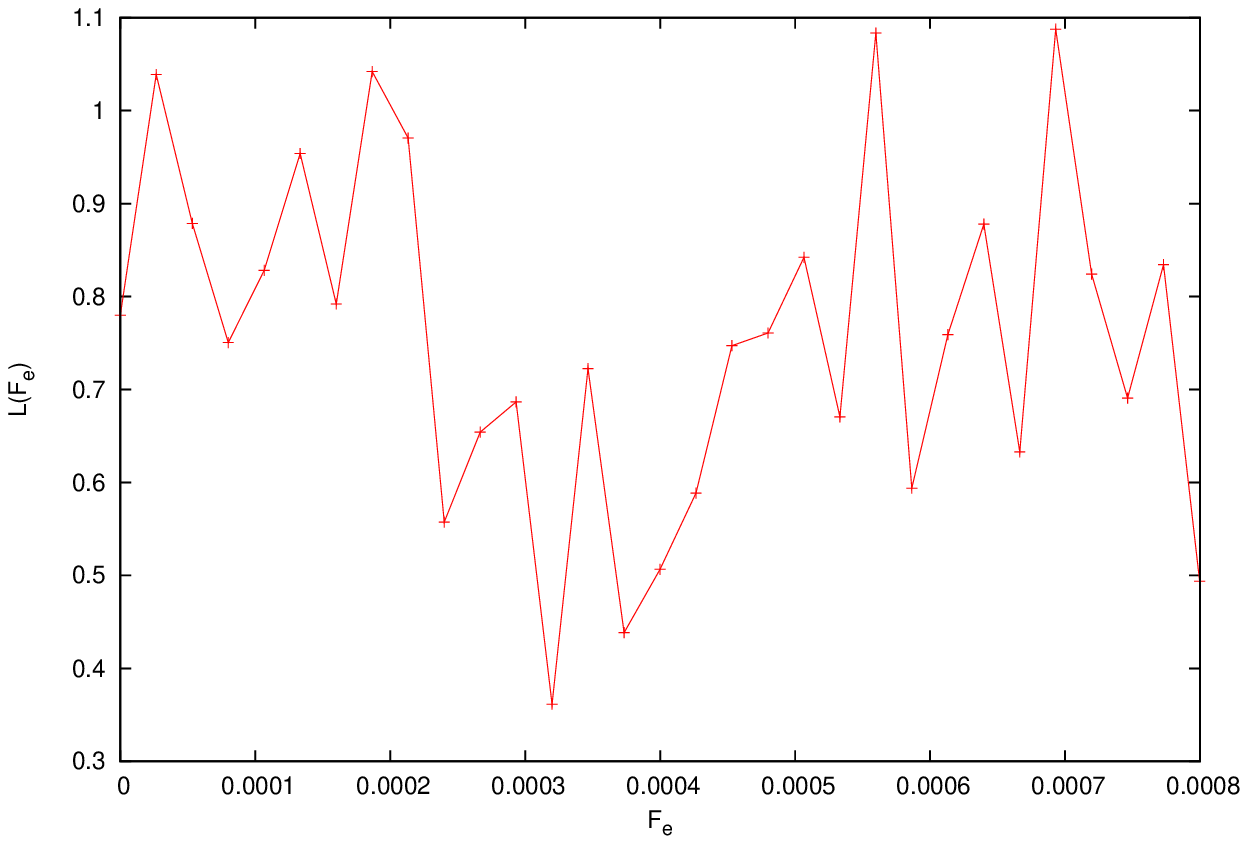}
\hspace{2mm}
\includegraphics[width=6cm]{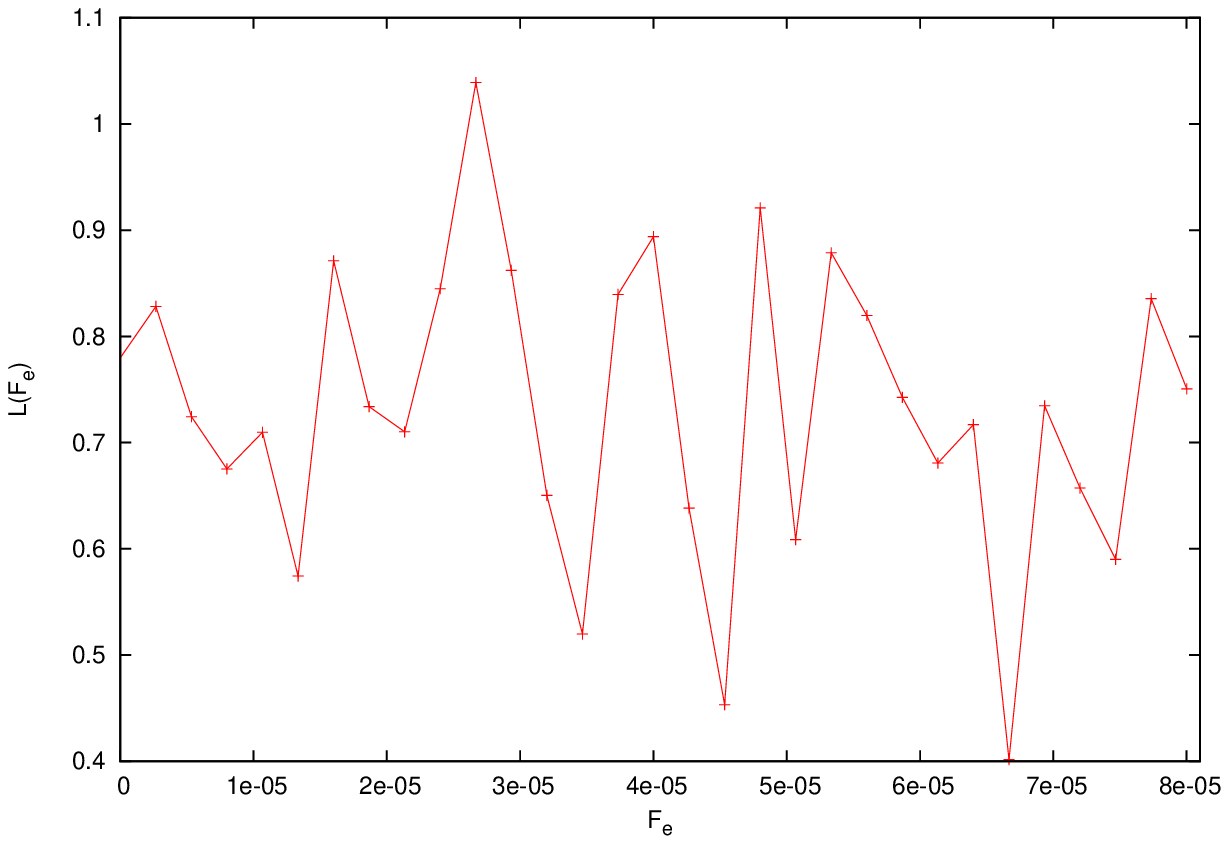}
\hspace{2mm}
\includegraphics[width=6cm]{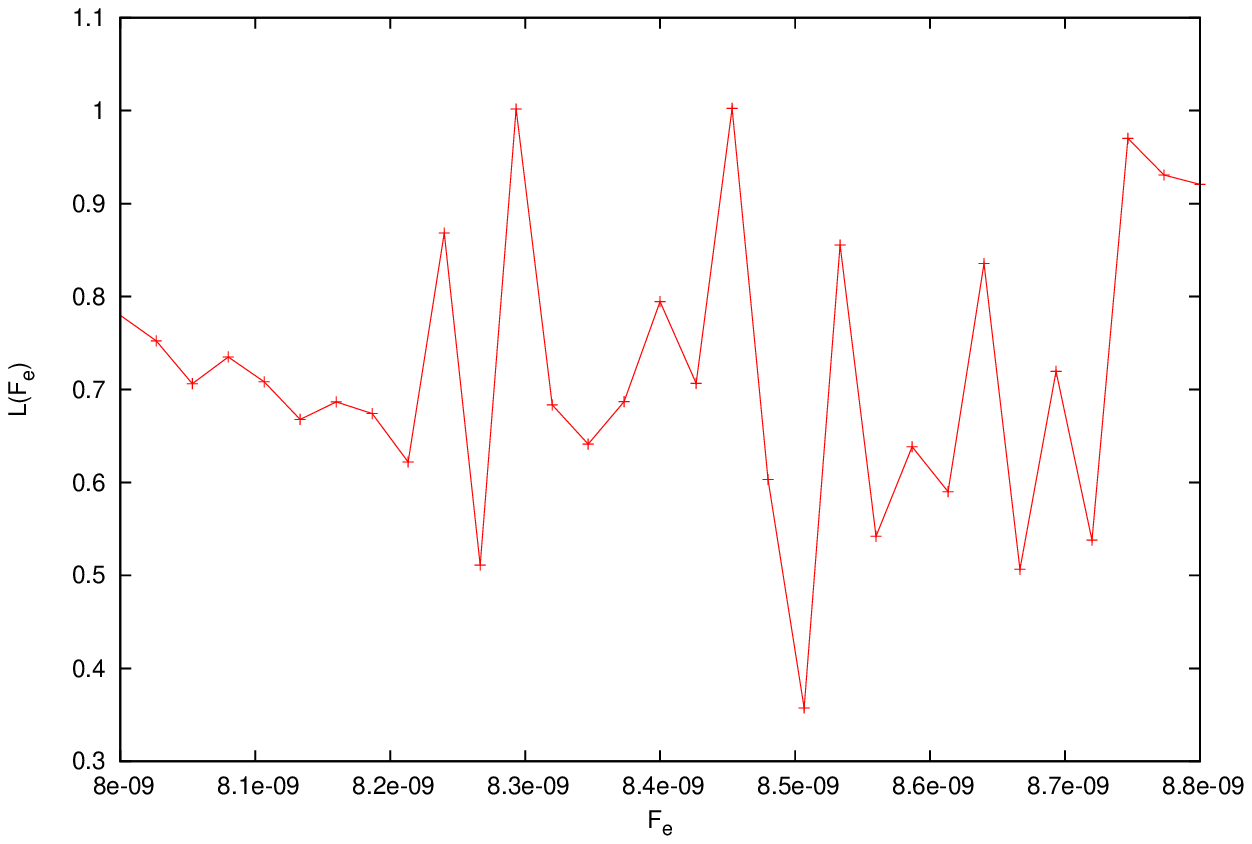}
\hspace{2mm}
\includegraphics[width=6cm]{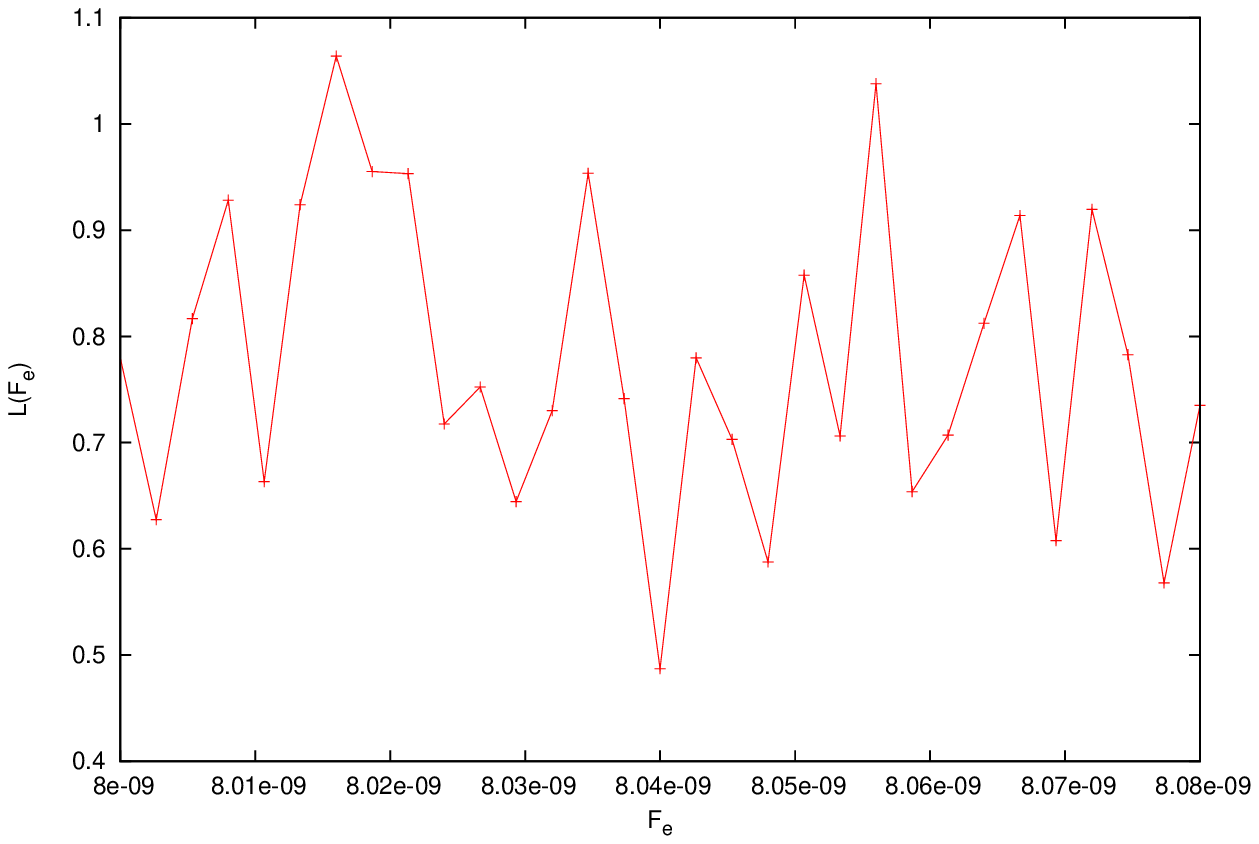}
\caption{Numerical computation of the coefficient $L(F_e)$ as a function of the bias $F_e$ at different orders of magnitude of the bias. The coefficient $L(F_e)$ presents a very irregular structure, which may prevent the existence of the limit $L(0)=\lim_{F_e \rightarrow 0} L(F_e)$.}\label{lintransp}
\end{figure}

\begin{figure}
\centering
\includegraphics[width=10.5cm]{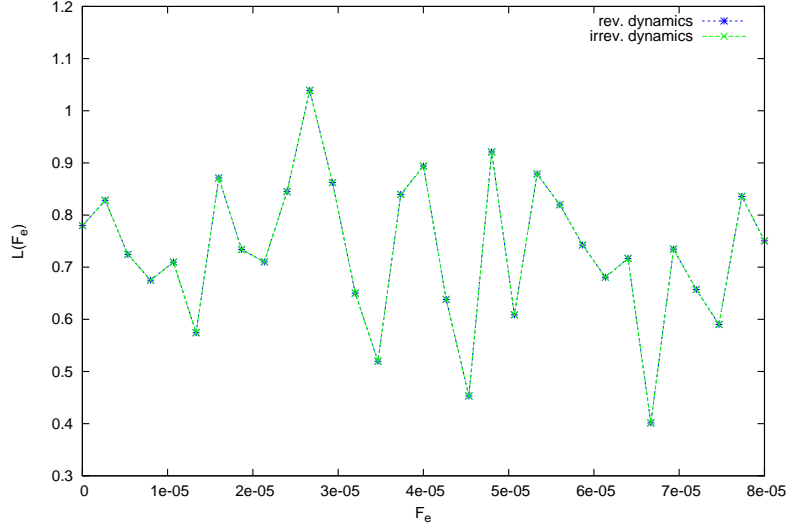}
\caption{Transport coefficient $L(F_e)$ for the irreversible, $K=M N$, and the corresponding reversible map, $K=M$, with $N_{ens}=5\cdot 10^5$. The irreversible component of the dynamics, $N$, leaves the response of the system unaffected.}\label{transport}
\end{figure}

To study the linear response for the irreversible map, that is the existence of the
limit $L(0)=\lim_{F_e \rightarrow 0} L(F_e)$,
we varied the value of bias within four different windows of magnitude, cf.\ Fig.\ref{lintransp},
and obtained that, in spite of the validity of the $\Lambda$-FR, no linear
response can be claimed for our low-dimensional system, unless this is verified at exceedingly small bias.
This fact cannot be blamed on the irreversible nature of the evolution, since we have observed
that the response of the irreversible map coincides with the response of the corresponding T-symmetric
one, cf.\ Fig.\ref{transport}. It is more related to the irregularity typical of transport phenomena
in low dimension \cite{RKbook}.
Therefore, the dynamics of the map $K=MN$ proves that its irreversible component,
the map $N$, affects neither the validity of the $\Lambda$-FR nor the response of the system to an
external bias.

{\em Per se},
the fact that neither the approach of \cite{Gall98} nor that of SRE are applicable does not imply
that no linear response can be established. However, this a clear observation for our model, which
does not enjoy many properties of the systems of \cite{Gall98} and many others of the systems of
\cite{SRE}. More importantly, as pointed out in e.g.\ Refs.\cite{JR,Ladek,RCPhysD2002}, the physical linear response
relies on the occurrence of \textit{local equilibrium}, in the sense that real space may be thought of as a ``collection" of cells, each of which contains a statistically significant number of interacting particles.
Clearly, our two-dimensional dynamical system may mimic only a few features of a real $N$-particle system,
and the local equilibrium property is out of question.

\section{Conclusions}
\label{sec:sec4}
In this paper, we have studied a kind of baker model whose properties are determined by two parameters,
one of which, $q$, may be conveniently tuned in order to fix the distance from the ``equilibrium'' state.
In particular, we have generalized the nonequilibrium map introduced in \cite{CKDR}, which can be
recovered from our present model by a suitable choice of $q$.
The system studied in this work is only a caricature of a real particle system subjected to
the action of an external driving, but can be studied in detail and thus is useful in understanding
the projection procedures typically thought be necessary to obtain a coarser, stochastic-like, description
from a microscopic, deterministic, one.

If one is interested in one specific observable, e.g. $\Lambda$, our model allows a simple
identification of the relevant and of the irrelevant variables: the phase function $\zL$ depends
on the Jacobians of the mapping, and the Jacobians only depend on one of the variables, $x$, which is then
the only relevant variable. This implies that adding a source of irreversibility which affects the
``irrelevant" degree of freedom $y$, has no effect on the equilibrium state defined by
$\langle \zL \rangle=0$, and on the validity of the $\Lambda$-FR.

Projecting the invariant measure on the reduced space of the relevant variable not only produces a
probability density, which is smooth along the unstable manifold, except for one point of discontinuity,
but also shows that equilibrium survives in the reduced space of the $x$-variable, even in presence
of irreversible full phase space dynamics. More precisely, equilibrium survives in the form of
detailed balance, which is the notion characterizing the equilibrium states in the spaces of relevant
observables, e.g.\ the $\mu$-space. This is due to the fact that in these spaces the details of the
dynamics of the ``noisy'' degrees of freedom are irrelevant.

We have also considered the validity of the $\zL$-FR and of the transport properties of an irreversible
dynamical system. So far the validity of the $\zL$-FR has been derived as a property of T-symmetric
dynamical systems. Our analysis, also supported by numerical tests, shows that as long as the source
of irreversibility affects only the ``irrelevant" degrees of freedom, the$\zL$-FR and the transport
laws (linear and nonlinear response) hold indistinguishably for both reversible and irreversible
phase space dynamics. Thus, our results extend to phase space dynamics some of the considerations
raised e.g.\ in Ref.\cite{Jona} for stochastic dynamics. This is, in fact, done by relating the
phase space dynamics to its projections which, according to basic tenets of statistical mechanics,
should result in stochastic evolutions.\\
In our investigation on a prototype of irreversible dynamical system, we discussed a model which is, manifestly, non Anosov: non-invertibility is merely one way to accomplish that. However, one may think of more realistic non-invertible dynamics \cite{JV}. In that case, the orientations of the phase space regions in which the phase space is decomposed will play an important role. Nevertheless, in more physical, higher dimensional, models, it might not be so straightforward to identify relevant or irrelevant variables with fixed directions in phase space, as was possible to do with our model. Investigations of further, more elaborate models is, therefore, worthwhile.

\vskip 10pt
\noindent
{\bf Acknowledgments}

\vskip 5pt
The Authors gratefully aknowledge fruitful discussions with Rainer Klages.

\end{document}